\documentclass[notitlepage,10pt]{article}
\usepackage{tikz,amssymb,mathrsfs,amsmath,empheq,esint,enumitem,graphicx,color,comment}

\usepackage{hyperref}
\usepackage[title]{appendix}

\catcode`\@=11
\@addtoreset{equation}{section}

\catcode`\@=12

\newtheorem{Theorem}{Theorem}[section]
\newtheorem{Lemma}[Theorem]{Lemma}
\newtheorem{Proposition}[Theorem]{Proposition}
\newtheorem{Corollary}[Theorem]{Corollary}
\newtheorem{Remark}[Theorem]{Remark}

\def\sech{\mathrm{\hspace{1.5pt}sech\hspace{1.5pt}}}
\def\R{{\mathbb R}}
\def\N{{\mathbb N}}
\def\S{{\mathbb S}}
\def\Z{{\mathbb Z}}
\def\f{{\varphi}}
\def\e{{\varepsilon}}
\def\veps{\varphi_{\varepsilon}}
\def\RT{{\mathcal{R}}}

\def\E{{\mathcal E}}
\def\F{{\mathcal F}}
\def\G{{\mathcal G}}
\def\L{{\mathcal L}}
\def\M{{\mathcal M}}
\def\T{{\mathcal{T}}}
\def\V{{\mathcal{V}}}
\def\X{{\mathcal{X}}}
\def\normal{{N}}
\def\am{{\mathrm{am}}}
\def\dn{{\mathrm{dn}}}

\def\qed{\hfill {$\square$}\goodbreak \medskip}

\title{Annular type surfaces with fixed boundary and\\ with prescribed, almost constant mean curvature}
\author{Paolo Caldiroli\footnote{\href{mailto:paolo.caldiroli@unito.it}{paolo.caldiroli@unito.it}}, Gabriele Cora\footnote{\href{mailto:gabriele.cora@unito.it}{gabriele.cora@unito.it}}, Alessandro Iacopetti\footnote{\href{mailto:alessandro.iacopetti@unito.it}{alessandro.iacopetti@unito.it}} \footnote{The authors are members of the Gruppo Nazionale per l'Analisi Matematica, la Probabilit\`a e le loro Applicazioni (GNAMPA) of the Istituto Nazionale di Alta Matematica (INdAM).}
\vspace{9pt}
\\
\normalsize{Dipartimento di Matematica ``Giuseppe Peano''}
\\
\normalsize{Universit\`a degli Studi di Torino}\\ 
\normalsize{via Carlo Alberto, 10 -- 10123 Torino, Italy}
}

\date{}

\begin{document}
\maketitle
\begin{abstract}
\noindent
We prove existence and nonexistence results for annular type parametric surfaces with prescribed, almost constant mean curvature, characterized as normal graphs of compact portions of unduloids or nodoids in $\mathbb{R}^{3}$, and whose boundary consists of two coaxial circles of the same radius. 
\smallskip

\smallskip

\end{abstract}

\section{Introduction}

In this work, we deal with the existence of immersed surfaces in \(\mathbb{R}^{3}\), of annular type, tethered to two parallel planes along two coaxial circles of the same radius and with prescribed mean curvature, close to a nonzero constant $H_{0}$ (henceforth we
assume $H_{0}=1$). 

In recent years, this issue has attracted quite an interest, especially in the case of constant mean curvature (see, e.g., \cite{KoisoMiyamoto, KoisoPalmerPiccione15, KoisoPalmerPiccione17, Patnaik94, Wente} and related references therein). In fact, constant mean curvature surfaces well approximate static equilibrium configurations of continua supported by surface tension, neglecting external forces. However, the presence of external forces, like gravity or other fields, may produce pressure differences which locally change, hence, according to the Young-Laplace equation, the mean curvature along phase interfaces is no longer a constant (see \cite{Finn, Fitzpatrick}). This motivates a study in the case of prescribed, near constant but not necessarily constant mean curvature. 

Actually, in this study, the main role is played by some classes of Delaunay surfaces, more precisely by unduloids and nodoids. Let us recall that Delaunay surfaces are complete, axially symmetric surfaces of constant mean curvature. They can be divided into five different types: planes, cathenoids, spheres, unduloids (including cylinders) and nodoids. We focus our attention on the last two families, made by noncompact surfaces with nonzero constant mean curvature. Unduloids and nodoids can be obtained by rotating the roulette of an ellipse, called \emph{elliptic catenary}, or the roulette of a hyperbola, called \emph{nodary}, respectively, about the revolution axis (see \cite{Delaunay, Eells, Kenmotsu} for more details). We also point out that unduloids and nodoids play a very important role in the construction of complete, constant mean curvature surfaces, with prescribed genus and possibly with ends, in an intriguing line of research originated from some breakthrough works by Kapouleas (see \cite{Kapouleas90, Kapouleas91}) and pursued by many other authors (see \cite{BreinerKapouleas, Grosse1993, GrosseKusnerSullivan2007, MazzeoPacard2001, MazzeoPacardPollack2001} and the bibliography therein).

Here we investigate the following problem: let us fix two parallel planes $\Pi_-$ and $\Pi_+$ in $\mathbb{R}^3$ and let us consider the families of sections of unduloids and nodoids with mean curvature one, whose axis of revolution is orthogonal to the planes $\Pi_-$ and $\Pi_+$, which are symmetric with respect to reflection around the plane half way between $\Pi_-$ and $\Pi_+$ and whose boundaries are couples of circles of the same radius lying on $\Pi_-$ and on $\Pi_+$. They form two real analytic 3-parameter families of compact surfaces which can be properly parameterized as follows. Let us introduce a Cartesian system such that $\Pi_\pm=\{(x,y,z)\in\mathbb{R}^3\mid z=\pm L\}$, for some $L>0$. Then, let us consider Delaunay roulettes defined as curves lying in the plane $y=0$, parameterized by the solution $(x_{a}(t),z_{a}(t))$ to the initial problem
\begin{equation}\label{eq:xaza}
\left\{\begin{array}{l}
x''=(1+2\gamma_a)x-2x^{3}\\
z'=x^{2}-\gamma_a
\end{array}\right.
\quad\left\{\begin{array}{l}
x(0)=1+a\,,\quad x'(0)=0\\
z(0)=0
\end{array}\right.
\end{equation}
where 
\[
\gamma_a:=a(1+a)
\]
and $a\in(-1,0)\cup(0,\infty)$ is a parameter, called Delaunay parameter. We point out that the function $x_a$ is even, positive and periodic of some period $2\tau_a$ (see \cite[Lemma 2.3]{CIM}). Such mapping can be interpreted as the speed of the parametrization, since the following identity holds:
\begin{equation}\label{eq:isothermal}
x_a=\sqrt{(x_a')^2+(z_a')^2}\,.
\end{equation}

When $a\in(-1,0)$ the curve is simple and is an elliptic catenary; the corresponding surface of revolution is embedded and is an unduloid. In particular, for $a=-\frac12$ such surface is a cylinder. Instead, when $a>0$, the curve has self-intersections and is a nodary; the corresponding surface of revolution is immersed but not embedded and is a nodoid. In both cases $\min\{|a|,1+a\}$ measures the neck-size, namely the distance of the surface from the axis of revolution. See the figure below. 
\bigskip


\noindent
\begin{tabular}{cc}
\!\!\begin{tikzpicture}[scale=1.1]
\draw[->](-1.2767,0)--(0.2+3*1.2767,0);
\draw[->](-1.2767,0)--(-1.2767,1.5);
\draw[dashed](-1.2767,0.1000)--(3*1.2767,0.1000);
\draw[dashed](-1.2767,0.9000)--(3*1.2767,0.9000);
\filldraw(-1.2767,0.1000)circle(0.4pt)node[left]{\footnotesize{$\min\{|a|,1+a\}$}};
\filldraw(-1.2767,0.9000)circle(0.4pt)node[left]{\footnotesize{$\max\{|a|,1+a\}$}};
\draw[thick,smooth,xshift=-1.2767cm](0.0000,0.9000)--(0.0280,0.8997)--(0.0429,0.8992)--(0.0728,0.8977)--(0.1025,0.8954)--(0.1320,0.8923)--(0.1612,0.8885)--(0.1756,0.8863)--(0.1900,0.8839)--(0.2184,0.8787)--(0.2324,0.8758)--(0.2463,0.8728)--(0.2738,0.8663)--(0.3007,0.8592)--(0.3270,0.8516)--(0.3527,0.8436)--(0.3778,0.8350)--(0.4022,0.8261)--(0.4259,0.8168)--(0.4490,0.8071)--(0.4713,0.7973)--(0.4930,0.7871)--(0.5140,0.7768)--(0.5343,0.7663)--(0.5539,0.7557)--(0.5821,0.7396)--(0.6088,0.7233)--(0.6341,0.7071)--(0.6580,0.6908)--(0.6879,0.6694)--(0.7090,0.6535)--(0.7353,0.6326)--(0.7597,0.6122)--(0.7823,0.5923)--(0.8085,0.5682)--(0.8370,0.5405)--(0.8667,0.5099)--(0.8932,0.4811)--(0.9233,0.4466)--(0.9525,0.4115)--(0.9781,0.3794)--(1.0008,0.3500)--(1.0232,0.3204)--(1.0434,0.2934)--(1.0663,0.2644)--(1.0870,0.2362)--(1.1089,0.2091)--(1.1258,0.1894)--(1.1439,0.1700)--(1.1620,0.1526)--(1.1775,0.1394)--(1.1876,0.1318)--(1.2039,0.1212)--(1.2193,0.1132)--(1.2351,0.1069)--(1.2497,0.1029)--(1.2629,0.1008)--(1.2712,0.1001)--(1.2767,0.1000);
\draw[thick,smooth,xshift=1.2767cm](-0.0000,0.9000)--(-0.0280,0.8997)--(-0.0429,0.8992)--(-0.0728,0.8977)--(-0.1025,0.8954)--(-0.1320,0.8923)--(-0.1612,0.8885)--(-0.1756,0.8863)--(-0.1900,0.8839)--(-0.2184,0.8787)--(-0.2324,0.8758)--(-0.2463,0.8728)--(-0.2738,0.8663)--(-0.3007,0.8592)--(-0.3270,0.8516)--(-0.3527,0.8436)--(-0.3778,0.8350)--(-0.4022,0.8261)--(-0.4259,0.8168)--(-0.4490,0.8071)--(-0.4713,0.7973)--(-0.4930,0.7871)--(-0.5140,0.7768)--(-0.5343,0.7663)--(-0.5539,0.7557)--(-0.5821,0.7396)--(-0.6088,0.7233)--(-0.6341,0.7071)--(-0.6580,0.6908)--(-0.6879,0.6694)--(-0.7090,0.6535)--(-0.7353,0.6326)--(-0.7597,0.6122)--(-0.7823,0.5923)--(-0.8085,0.5682)--(-0.8370,0.5405)--(-0.8667,0.5099)--(-0.8932,0.4811)--(-0.9233,0.4466)--(-0.9525,0.4115)--(-0.9781,0.3794)--(-1.0008,0.3500)--(-1.0232,0.3204)--(-1.0434,0.2934)--(-1.0663,0.2644)--(-1.0870,0.2362)--(-1.1089,0.2091)--(-1.1258,0.1894)--(-1.1439,0.1700)--(-1.1620,0.1526)--(-1.1775,0.1394)--(-1.1876,0.1318)--(-1.2039,0.1212)--(-1.2193,0.1132)--(-1.2351,0.1069)--(-1.2497,0.1029)--(-1.2629,0.1008)--(-1.2712,0.1001)--(-1.2767,0.1000);
\draw[thick,smooth,xshift=1.2767cm](0.0000,0.9000)--(0.0280,0.8997)--(0.0429,0.8992)--(0.0728,0.8977)--(0.1025,0.8954)--(0.1320,0.8923)--(0.1612,0.8885)--(0.1756,0.8863)--(0.1900,0.8839)--(0.2184,0.8787)--(0.2324,0.8758)--(0.2463,0.8728)--(0.2738,0.8663)--(0.3007,0.8592)--(0.3270,0.8516)--(0.3527,0.8436)--(0.3778,0.8350)--(0.4022,0.8261)--(0.4259,0.8168)--(0.4490,0.8071)--(0.4713,0.7973)--(0.4930,0.7871)--(0.5140,0.7768)--(0.5343,0.7663)--(0.5539,0.7557)--(0.5821,0.7396)--(0.6088,0.7233)--(0.6341,0.7071)--(0.6580,0.6908)--(0.6879,0.6694)--(0.7090,0.6535)--(0.7353,0.6326)--(0.7597,0.6122)--(0.7823,0.5923)--(0.8085,0.5682)--(0.8370,0.5405)--(0.8667,0.5099)--(0.8932,0.4811)--(0.9233,0.4466)--(0.9525,0.4115)--(0.9781,0.3794)--(1.0008,0.3500)--(1.0232,0.3204)--(1.0434,0.2934)--(1.0663,0.2644)--(1.0870,0.2362)--(1.1089,0.2091)--(1.1258,0.1894)--(1.1439,0.1700)--(1.1620,0.1526)--(1.1775,0.1394)--(1.1876,0.1318)--(1.2039,0.1212)--(1.2193,0.1132)--(1.2351,0.1069)--(1.2497,0.1029)--(1.2629,0.1008)--(1.2712,0.1001)--(1.2767,0.1000);
\draw[thick,smooth,xshift=3*1.2767cm](-0.0000,0.9000)--(-0.0280,0.8997)--(-0.0429,0.8992)--(-0.0728,0.8977)--(-0.1025,0.8954)--(-0.1320,0.8923)--(-0.1612,0.8885)--(-0.1756,0.8863)--(-0.1900,0.8839)--(-0.2184,0.8787)--(-0.2324,0.8758)--(-0.2463,0.8728)--(-0.2738,0.8663)--(-0.3007,0.8592)--(-0.3270,0.8516)--(-0.3527,0.8436)--(-0.3778,0.8350)--(-0.4022,0.8261)--(-0.4259,0.8168)--(-0.4490,0.8071)--(-0.4713,0.7973)--(-0.4930,0.7871)--(-0.5140,0.7768)--(-0.5343,0.7663)--(-0.5539,0.7557)--(-0.5821,0.7396)--(-0.6088,0.7233)--(-0.6341,0.7071)--(-0.6580,0.6908)--(-0.6879,0.6694)--(-0.7090,0.6535)--(-0.7353,0.6326)--(-0.7597,0.6122)--(-0.7823,0.5923)--(-0.8085,0.5682)--(-0.8370,0.5405)--(-0.8667,0.5099)--(-0.8932,0.4811)--(-0.9233,0.4466)--(-0.9525,0.4115)--(-0.9781,0.3794)--(-1.0008,0.3500)--(-1.0232,0.3204)--(-1.0434,0.2934)--(-1.0663,0.2644)--(-1.0870,0.2362)--(-1.1089,0.2091)--(-1.1258,0.1894)--(-1.1439,0.1700)--(-1.1620,0.1526)--(-1.1775,0.1394)--(-1.1876,0.1318)--(-1.2039,0.1212)--(-1.2193,0.1132)--(-1.2351,0.1069)--(-1.2497,0.1029)--(-1.2629,0.1008)--(-1.2712,0.1001)--(-1.2767,0.1000);
\end{tikzpicture}
&
\!\!\begin{tikzpicture}[scale=1.1]
\draw[->](-2*0.6403,0)--(0.2+6*0.6403,0);
\draw[white](-2*0.6403,0)--(-2*0.6403,-0.135);
\draw[->](-2*0.6403,0)--(-2*0.6403,1.55);
\draw[dashed](-2*0.6403,0.1727)--(6*0.6403,0.1727);
\draw[dashed](-2*0.6403,1.1727)--(6*0.6403,1.1727);
\filldraw(-2*0.6403,0.1727)circle(0.4pt)node[left]{\footnotesize{$a$}};
\filldraw(-2*0.6403,1.1727)circle(0.4pt)node[left]{\footnotesize{$1$+$a$}};
\draw[thick,smooth](0,1.1727)--(0.0650,1.1703)--(0.1292,1.1630)--(0.1918,1.1513)--(0.2402,1.1389)--(0.2983,1.1199)--(0.3531,1.0977)--(0.4043,1.0728)--(0.4519,1.0457)--(0.5119,1.0051)--(0.5646,0.9626)--(0.6103,0.9194)--(0.6443,0.8824)--(0.6739,0.8459)--(0.7111,0.7926)--(0.7436,0.7364)--(0.7725,0.6735)--(0.7954,0.6066)--(0.8125,0.5264)--(0.8178,0.4673)--(0.8155,0.3998)--(0.7985,0.3196)--(0.7779,0.2713)--(0.7454,0.2250)--(0.7047,0.1918)--(0.6632,0.1755)--(0.6403,0.1727);
\draw[thick,smooth](0,1.1727)--(-0.0650,1.1703)--(-0.1292,1.1630)--(-0.1918,1.1513)--(-0.2402,1.1389)--(-0.2983,1.1199)--(-0.3531,1.0977)--(-0.4043,1.0728)--(-0.4519,1.0457)--(-0.5119,1.0051)--(-0.5646,0.9626)--(-0.6103,0.9194)--(-0.6443,0.8824)--(-0.6739,0.8459)--(-0.7111,0.7926)--(-0.7436,0.7364)--(-0.7725,0.6735)--(-0.7954,0.6066)--(-0.8125,0.5264)--(-0.8178,0.4673)--(-0.8155,0.3998)--(-0.7985,0.3196)--(-0.7779,0.2713)--(-0.7454,0.2250)--(-0.7047,0.1918)--(-0.6632,0.1755)--(-0.6403,0.1727);
\draw[thick,smooth,xshift=-2*0.6403cm](0,1.1727)--(0.0650,1.1703)--(0.1292,1.1630)--(0.1918,1.1513)--(0.2402,1.1389)--(0.2983,1.1199)--(0.3531,1.0977)--(0.4043,1.0728)--(0.4519,1.0457)--(0.5119,1.0051)--(0.5646,0.9626)--(0.6103,0.9194)--(0.6443,0.8824)--(0.6739,0.8459)--(0.7111,0.7926)--(0.7436,0.7364)--(0.7725,0.6735)--(0.7954,0.6066)--(0.8125,0.5264)--(0.8178,0.4673)--(0.8155,0.3998)--(0.7985,0.3196)--(0.7779,0.2713)--(0.7454,0.2250)--(0.7047,0.1918)--(0.6632,0.1755)--(0.6403,0.1727);
\draw[thick,xshift=2*0.6403cm](0,1.1727)--(-0.0650,1.1703)--(-0.1292,1.1630)--(-0.1918,1.1513)--(-0.2402,1.1389)--(-0.2983,1.1199)--(-0.3531,1.0977)--(-0.4043,1.0728)--(-0.4519,1.0457)--(-0.5119,1.0051)--(-0.5646,0.9626)--(-0.6103,0.9194)--(-0.6443,0.8824)--(-0.6739,0.8459)--(-0.7111,0.7926)--(-0.7436,0.7364)--(-0.7725,0.6735)--(-0.7954,0.6066)--(-0.8125,0.5264)--(-0.8178,0.4673)--(-0.8155,0.3998)--(-0.7985,0.3196)--(-0.7779,0.2713)--(-0.7454,0.2250)--(-0.7047,0.1918)--(-0.6632,0.1755)--(-0.6403,0.1727);
\draw[thick,smooth,xshift=2*0.6403cm](0,1.1727)--(0.0650,1.1703)--(0.1292,1.1630)--(0.1918,1.1513)--(0.2402,1.1389)--(0.2983,1.1199)--(0.3531,1.0977)--(0.4043,1.0728)--(0.4519,1.0457)--(0.5119,1.0051)--(0.5646,0.9626)--(0.6103,0.9194)--(0.6443,0.8824)--(0.6739,0.8459)--(0.7111,0.7926)--(0.7436,0.7364)--(0.7725,0.6735)--(0.7954,0.6066)--(0.8125,0.5264)--(0.8178,0.4673)--(0.8155,0.3998)--(0.7985,0.3196)--(0.7779,0.2713)--(0.7454,0.2250)--(0.7047,0.1918)--(0.6632,0.1755)--(0.6403,0.1727);
\draw[thick,xshift=4*0.6403cm](0,1.1727)--(-0.0650,1.1703)--(-0.1292,1.1630)--(-0.1918,1.1513)--(-0.2402,1.1389)--(-0.2983,1.1199)--(-0.3531,1.0977)--(-0.4043,1.0728)--(-0.4519,1.0457)--(-0.5119,1.0051)--(-0.5646,0.9626)--(-0.6103,0.9194)--(-0.6443,0.8824)--(-0.6739,0.8459)--(-0.7111,0.7926)--(-0.7436,0.7364)--(-0.7725,0.6735)--(-0.7954,0.6066)--(-0.8125,0.5264)--(-0.8178,0.4673)--(-0.8155,0.3998)--(-0.7985,0.3196)--(-0.7779,0.2713)--(-0.7454,0.2250)--(-0.7047,0.1918)--(-0.6632,0.1755)--(-0.6403,0.1727);
\draw[thick,smooth,xshift=4*0.6403cm](0,1.1727)--(0.0650,1.1703)--(0.1292,1.1630)--(0.1918,1.1513)--(0.2402,1.1389)--(0.2983,1.1199)--(0.3531,1.0977)--(0.4043,1.0728)--(0.4519,1.0457)--(0.5119,1.0051)--(0.5646,0.9626)--(0.6103,0.9194)--(0.6443,0.8824)--(0.6739,0.8459)--(0.7111,0.7926)--(0.7436,0.7364)--(0.7725,0.6735)--(0.7954,0.6066)--(0.8125,0.5264)--(0.8178,0.4673)--(0.8155,0.3998)--(0.7985,0.3196)--(0.7779,0.2713)--(0.7454,0.2250)--(0.7047,0.1918)--(0.6632,0.1755)--(0.6403,0.1727);
\draw[thick,xshift=6*0.6403cm](0,1.1727)--(-0.0650,1.1703)--(-0.1292,1.1630)--(-0.1918,1.1513)--(-0.2402,1.1389)--(-0.2983,1.1199)--(-0.3531,1.0977)--(-0.4043,1.0728)--(-0.4519,1.0457)--(-0.5119,1.0051)--(-0.5646,0.9626)--(-0.6103,0.9194)--(-0.6443,0.8824)--(-0.6739,0.8459)--(-0.7111,0.7926)--(-0.7436,0.7364)--(-0.7725,0.6735)--(-0.7954,0.6066)--(-0.8125,0.5264)--(-0.8178,0.4673)--(-0.8155,0.3998)--(-0.7985,0.3196)--(-0.7779,0.2713)--(-0.7454,0.2250)--(-0.7047,0.1918)--(-0.6632,0.1755)--(-0.6403,0.1727);
\end{tikzpicture}
\\
\small{Fig.~\!1-a. Elliptic catenary ($-1<a<0$) }~&~~
\small{Fig.~\!1-b. Nodary ($a>0$)}
\end{tabular}
\bigskip

\noindent
Fixing $(p,q)\in\mathbb{R}^2$, the symmetric unduloid or nodoid with Delaunay parameter $a$ and whose revolution axis is the vertical line passing through $(p,q,0)$ is given by the image of $X_{a}+p\mathbf{e}_{1}+q\mathbf{e}_{2}$ where $\mathbf{e}_{1}$ and $\mathbf{e}_{2}$ are the first two vectors of the canonical basis in $\mathbb{R}^{3}$ and 
\begin{equation}\label{Xa}
X_{a}(t,\theta)=\left[\begin{array}{c}x_a(t)\cos\theta\\x_a(t)\sin\theta\\z_a(t)\end{array}\right],\quad(t,\theta)\in\mathbb{R}\times[-\pi,\pi]\,.
\end{equation}
The triple $(a,p,q)\in((-1,0)\cup(0,\infty))\times\mathbb{R}\times\mathbb{R}$ parameterizes, with analytic dependence, the elements of the family of symmetric vertical unduloids, when $a\in(-1,0)$, and those ones of the family of vertical nodoids, when $a \in (0, \infty)$. In view of \eqref{eq:xaza}, for every $(a,p,q)\in((-1,0)\cup(0,\infty))\times\mathbb{R}\times\mathbb{R}$, the mean curvature of $X_{a}+p\mathbf{e}_{1}+q\mathbf{e}_{2}$ at every point is 1. 

Now we fix $a\in (-1,0)\cup(0,\infty)$, we take $T>0$ such that $z_a(T)=L$, and we consider the family of compact surfaces $\Sigma_{a,T,p,q}$ parameterized by $X_{a}(t,\theta)+p\mathbf{e}_{1}+q\mathbf{e}_{2}$ with $(t,\theta)\in[-T,T]\times[-\pi,\pi]$. The boundary of $\Sigma_{a,T,p,q}$ consists of the couple of circles lying on the planes $z=\pm z_a(T)$, centered at $(p,q,\pm z_a(T))$, and with radius $x_a(T)$. 
Then we consider surfaces (of annular type) which are \emph{normal graphs of $\Sigma_{a,T,p,q}$ with fixed boundary}. More explicitly, we introduce the Gauss map associated to $X_{a}$, defined by
\[
N_{a}(t,\theta):=\frac{\partial_t X_{a}(t,\theta)\wedge\partial_\theta X_{a}(t,\theta)}{|\partial_t X_{a}(t,\theta)\wedge\partial_\theta X_{a}(t,\theta)|}
\]
and we consider surfaces with parametrization of the form
\begin{equation}
\label{normal-graph}
X(t,\theta)=X_{a}(t,\theta)+p\mathbf{e}_{1}+q\mathbf{e}_{2}+\varphi(t,\theta)N_{a}(t,\theta)\,,\quad(t,\theta)\in[-T,T]\times[-\pi,\pi]
\end{equation}
with $\varphi\colon[-T,T]\times\mathbb{S}^{1}\to\mathbb{R}$ regular enough and such that $\varphi(\pm T,\cdot)=0$ on $\S^1$, with the identification $\mathbb{S}^{1}=\mathbb{R}/2\pi\mathbb{Z}$. 

Fixing $a\in(-1,0)\cup(0,\infty)$, $T>0$, and $(p,q)\in\mathbb{R}^{2}$, we investigate the stability of compact sections of unduloids and nodoids $\Sigma_{a,T,p,q}$ when we perturb the prescribed mean curvature function and we hold the boundary of the surfaces $\Sigma_{a,T,p,q}$. This problem can be phrased analytically as follows. Take a family of regular mappings $H_\varepsilon\colon\mathbb{R}^3\to\mathbb{R}$ depending smoothly on the parameter $\varepsilon\in\mathbb{R}$ and such that $H_0\equiv 1$. Look for a (continuous) path of mappings $\varepsilon\mapsto\varphi_{\varepsilon}$ with $\varepsilon\in(-\bar\varepsilon,\bar\varepsilon)$, such that $\varphi_{0}=0$ and $\varphi_{\varepsilon}\colon[-T,T]\times\mathbb{S}^{1}\to\mathbb{R}$ solving the Dirichlet problem
\begin{equation}
\label{M=H}
\left\{\begin{array}{l}
\mathfrak{M}(X_{a}+p\mathbf{e}_{1}+q\mathbf{e}_{2}+\varphi_{\varepsilon}N_{a})=H_{\varepsilon}(X_{a}+p\mathbf{e}_{1}+q\mathbf{e}_{2}+\varphi_{\varepsilon}N_{a})\text{ in }[-T,T]\times\mathbb{S}^{1}\\
\varphi_{\varepsilon}(\pm T,\cdot)=0
\end{array}
\right.
\end{equation}
for every $\varepsilon\in(-\bar\varepsilon,\bar\varepsilon)$, where 
$\mathfrak{M}$ is the mean curvature operator (defined in Section \ref{S:preliminaries}). 

We will assume that the mappings $H_{\varepsilon}$ are even with respect to the last variable. This hypothesis turns out to be technically very useful since it helps to reduce degeneracy (see below) and it allows us to set up problem \eqref{M=H} in functional spaces of mappings $\varphi(t,\theta)$ which are even with respect to $t$. This means that the surfaces parametrized by \eqref{normal-graph} are symmetric with respect to the horizontal plane $z=0$.  

For future reference, it is convenient to introduce here the functional spaces we will be dealing with, defined in this way:
\[
\mathcal{X}_{T}:=\{\varphi\in C^{2,\alpha}([-T,T]\times\mathbb{S}^{1},\mathbb{R})\mid \varphi(\pm T,\cdot)=0\,\text{ on }\S^1,\ \varphi(t,\cdot)=\varphi(-t,\cdot)\ \forall\, t\in[-T,T]\}
\]
(with $\alpha\in(0,1)$ fixed), and equipped with the standard $C^{2,\alpha}$-H\"older norm.

Before stating our results, we introduce the concept of degeneracy/nondegeneracy for compact sections of unduloids and nodoids. To this aim, we write down the linearized problem associated to \eqref{M=H}, given by
\begin{equation}
\label{L=0}
\left\{\begin{array}{l}
\mathfrak{L}_{a}\varphi=0\text{ in }[-T,T]\times\mathbb{S}^{1}\\
\varphi\in\mathcal{X}_{T}
\end{array}
\right.
\end{equation}
where $\mathfrak{L}_{a}$ is the Jacobi operator associated to the Delaunay surface $X_{a}$, defined by
\[
\mathfrak{L}_{a}\varphi:=\lim_{s\to 0}\frac{\mathfrak{M}(X_{a}+s\varphi N_{a})-\mathfrak{M}(X_{a})}{s}\,.
\]
We say that a surface $\Sigma_{a,T,p,q}$ is \emph{degenerate} (\emph{nondegenerate}, respectively) if 
$\mathfrak{L}_{a}$ has a nontrivial kernel (a null kernel, respectively) in $\mathcal{X}_{T}$. As we will see, degeneracy or not of $\Sigma_{a,T,p,q}$ strongly affects the question of existence of solutions to \eqref{M=H}. In fact, this property is independent of $(p,q)\in\mathbb{R}^{2}$, whereas, for  fixed $a$, the value of $T$ plays a fundamental role. For this reason, first of all, it is important to study the set
\[
\mathcal{T}_{a}:=\{T>0\mid\text{\eqref{L=0} admits a nontrivial solution}\}\,.
\]
More precisely, it is known that for every integer $j\ge 0$ the equation $\mathfrak{L}_{a}\varphi=0$ on $\mathbb{R}\times\mathbb{S}^{1}$ admits non null solutions of the form
\begin{equation}
\label{w}
\varphi(t,\theta)=w(t)\cos(j\theta)\quad\text{with $j\in\N$}
\end{equation}
where $w$ is an even solution of a certain ODE depending also on the parameter $a$ (see \eqref{eq:j-eq}).

We need to understand the behavior of this kind of solutions and, in particular, the structure of the nodal set of $w(t)$, because zeros of $w(t)$, if any, are elements of $\T_{a}$ and correspond to the degeneracy case. In order to achieve this goal, we use also the Floquet theory, and we attain the following result. 

\begin{Theorem}
\label{T:Ta} 
For every $a\in(-1,0)\cup(0,\infty)$ and for every $j\in\N$ let
\begin{equation}
\label{Taj}
\mathcal{T}_{a,j}:=\{T>0\mid \text{$\exists\, \varphi\in C^{2}(\mathbb{R}\times\mathbb{S}^{1})\setminus\{0\}$ as in \eqref{w} s.t. $ \mathfrak{L}_{a}\varphi=0$ on $\mathbb{R}\times\mathbb{S}^{1}$ and $\varphi( T,\cdot)=0$}\}\,.
\end{equation}
Then $\mathcal{T}_{a}=\bigcup_{j=0}^{\infty}\mathcal{T}_{a,j}$. Moreover:
\begin{itemize}
\item[$(i)$]
$\mathcal{T}_{a,0}=\{T_{a,k}\}_{k=0}^{\infty}$ with $T_{a,k}\in\big(k\tau_{a},\big(k+\frac12\big)\tau_{a}\big)$ for every $k\in \N$.
\item[$(ii)$]
If $a<0$, then $\mathcal{T}_{a,j}=\varnothing$ for every $j\ge 1$ and $\mathcal{T}_{a}=\mathcal{T}_{a,0}$. 
\item[$(iii)$]
If $a>0$, then $\mathcal{T}_{a,1}=\big\{(k+\frac12)\tau_{a}\big\}_{k=0}^{\infty}$ and $\mathcal{T}_{a,0}\cap\mathcal{T}_{a,1}=\varnothing$.
\item[$(iv)$]
If $a\in\big(0,\frac12(\sqrt3-1)\big]$, then $\mathcal{T}_{a,j}=\varnothing$ for every $j\ge 2$ {and $\mathcal{T}_{a}=\mathcal{T}_{a,0}\cup \T_{a,1}$}.
\item[$(v)$] for every $j\geq2$ there exists $a^*_j>0$ such that $ \T_{a,j} \neq \emptyset$ for $a \geq a^*_j$. 
\end{itemize}
\end{Theorem}

In the description of the set $\mathcal{T}_{a}$ stated by Theorem \ref{T:Ta} the most important subsets are $\mathcal{T}_{a,0}$ and $\mathcal{T}_{a,1}$, since they have a quite clear geometrical meaning. Indeed, if $T\not\in \mathcal{T}_{a,0}$, then, considering the class of surfaces $\Sigma_{a',T',p,q}$ with fixed $(p,q)\in\mathbb{R}^{2}$, and $(a',T')$ in a neighborhood of $(a,T)$, only the surface $\Sigma_{a,T,p,q}$ in this class has boundary consisting of the circles of radius $x_{a}(T)$, placed at a distance $2z_{a}(T)$ (see Remark \ref{R:Ta0}). In other words, when $T\in \mathcal{T}_{a,0}$, the surface $\Sigma_{a,T,p,q}$ is somehow less robust with respect to perturbations. 

Regarding the set $\mathcal{T}_{a,1}$, when $a>0$ and $T\in\mathcal{T}_{a,1}$, the surfaces $\Sigma_{a,T,p,q}$ are sections of nodoids which at the boundary intersect the planes $\Pi_{-}$ and $\Pi_{+}$ tangentially. A tangential contact with the end planes may cause a loss of stability. In fact, as we will see, a perturbation of the mean curvature has some effect. Note that tangent intersections cannot occur with unduloids, and in fact $\mathcal{T}_{a,1}$ is empty when $a<0$. 
\color{black}

Now, we go back to problem \eqref{M=H} and we exhibit our main results about it. When $T\not\in\mathcal{T}_{a}$ (nondegeneracy case) we can study \eqref{M=H} with the aid of the implicit function theorem and we obtain:

\begin{Theorem}
\label{T:existence}
Let $a\in(-1,0)\cup(0,\infty)$ and let $(\varepsilon,x,y,z)\mapsto H_{\varepsilon}(x,y,z)$ be a mapping in $C^{1}(\mathbb{R}\times\mathbb{R}^{3})$ such that:
\[
H_{0}(x,y,z)=1\quad\forall(x,y,z)\in\mathbb{R}^{3}\,,
\leqno{(H_{1})}
\]\vspace{-25pt}
\[
H_{\varepsilon}(x,y,z)=H_{\varepsilon}(x,y,-z)\quad\forall(x,y,z)\in\mathbb{R}^{3},~\forall\varepsilon\in\mathbb{R}\,.
\leqno{(H_{2})}
\]
\begin{itemize}
\item[$(i)$] If $T\in(0,\infty)\setminus \mathcal{T}_{a}$, then for every $(p,q)\in\mathbb{R}^{2}$ there exist $\bar{\varepsilon}=\bar{\varepsilon}(a,T,p,q)>0$ and a mapping $\varepsilon\mapsto\varphi_{\varepsilon}$ in $C^{1}((-\bar{\varepsilon},\bar{\varepsilon}),\mathcal{X}_{T})$ with $\varphi_{0}=0$ and $\varphi_{\varepsilon}$ solving \eqref{M=H} for every $\varepsilon\in(-\bar\varepsilon,\bar\varepsilon)$. 

\item[$(ii)$]If, in addition, $H_{\varepsilon}$ is independent of $x$ and $y$, then the same result holds for every $T\in(0,\infty)\setminus \mathcal{T}_{a,0}$ with $\varphi_{\varepsilon}$ depending just on $t$. 
\end{itemize}
\end{Theorem}

When $T\in\mathcal{T}_{a}$ (degeneracy case), the implicit function theorem does not work anymore because the operator $\mathfrak{L}_{a}$ has a non null kernel in $\mathcal{X}_{T}$. In this case different kind of results hold. When $T\in \mathcal{T}_{a,0}$, in general, no solution of \eqref{M=H} exists. More precisely, recalling that $\mathcal{T}_{a,0}=\{T_{a,k}\}_{k=0}^{\infty}$, we have:

\begin{Theorem}
\label{T:nonexistence}
Let $a\in(-1,0)\cup(0,\infty)$, $(p_0, q_0) \in \R^2$ and let $(\e, x, y, z) \mapsto H_{\varepsilon}(x,y,z)$ be a mapping in $C^{1}(\mathbb{R}\times\mathbb{R}^{3})$ satisfying $(H_{1})$. 
\begin{itemize}
\item[$(i)$]
If $T = T_{a,0}$ and 
\begin{equation}\label{eq:deHzero}
\left.\frac{\partial H_\varepsilon}{\partial \varepsilon}(x,y,z)\right|_{\varepsilon=0}\ne 0\quad\forall (x,y,z)\in\mathbb{R}^{3},
\end{equation}
then there cannot exist $\bar \e >0$, a continuous mapping $\varepsilon\mapsto(p_{\varepsilon},q_{\varepsilon}) \in \R^2$ and a $C^1$-mapping $\varepsilon\mapsto\varphi_{\varepsilon} \in  C^2([-T_{a,0},T_{a,0}]\times\S^1)$ with $|\e|\leq \bar \e$, such that $\left.(p_\e, q_\e)\right|_{\e = 0} = (p_0, q_0)$, $\varphi_{0}=0$ and for every $\varepsilon\in(-\bar\varepsilon,\bar\varepsilon)$ the triple $(p_\e, q_\e, \varphi_{\varepsilon})$ solves \eqref{M=H}. 
\item[$(ii)$] 
If $ T = T_{a,k}$, $k\geq1$, and 
\[
\int_{[-T_{a,k}, T_{a,k}]\times[-\pi, \pi]}x_{a}^{2}\overline{w}_{a,0}\,\left.\frac{\partial H_\varepsilon}{\partial \varepsilon}(X_a + p_0 \mathbf{e}_1 +q_0 \mathbf{e}_2 )\right|_{\varepsilon=0}\,dt\,d\theta\ne 0\,,
\]
where $\overline{w}_{a,0}$ is a nontrivial even solution of \eqref{L=0} depending just on $t$, then the same conclusion as in $(i)$ holds true. 
\end{itemize}
\end{Theorem}
In particular \eqref{eq:deHzero} holds true for $H_\e(x, y, z) = 1 + \e \widetilde H(x, y, z)$, with $\widetilde H \in C^1(\R^3)$ never vanishing. 

We notice that nonexistence is proved without asking for any symmetry condition on $H_{\varepsilon}$ or on the solutions of \eqref{M=H} with respect to $t$. 
To show Theorem \ref{T:nonexistence} we use the information stated in part $(i)$ of Theorem \ref{T:Ta}. Moreover, we consider a non null solution $\overline{w}_{a,0}$ of $\mathfrak{L}_{a}\varphi=0$ on $\mathbb{R}\times\mathbb{S}^{1}$ of the form \eqref{w} with $j=0$. Testing the equation in \eqref{M=H} with $\overline{w}_{a,0}$, writing down its expansion with respect to the smallness parameter $\varepsilon$, and analyzing carefully the different contributions coming from the right hand side and the left one, we obtain the result. 

We point out that for $a\in(-1,0)$ Theorems \ref{T:existence} and \ref{T:nonexistence} essentially cover all the cases, because, by Theorem \ref{T:Ta} part $(ii)$, either $T\not\in\mathcal{T}_{a}$ or $T\in\mathcal{T}_{a,0}$. 
When $a\in(0, \infty)$ the situation is different. In the degenerate case ($T\in\mathcal{T}_{a}$), besides when $T\in\mathcal{T}_{a,0}$, we can obtain some results only for $T\in 
\mathcal{T}_{a,1}$, provided that $a$ is not too large. In this case we have:

\begin{Theorem}
\label{T:existence2}
Let $a \in (0, \infty)$, $T \in \T_{a,1}$, $(p_0, q_0) \in \R^2$, and  let $H_{\varepsilon}(x,y,z)$ be a mapping in $C^{1}(\mathbb{R}\times\mathbb{R}^{3})$ satisfying $(H_{1})$.  
Let
\[
\widetilde H(x, y, z) := \left.\frac{\partial H_\e}{\partial \e}( x, y, z)\right|_{\e = 0}\,,
\]
and define 
\begin{equation}\label{eq:Q}
Q(x,y,z):=\frac12\left(\int_{0}^{x}\widetilde{H}(s,y,z)\,ds\,,\int_{0}^{y}\widetilde{H}(x,s,z)\,ds\,,0\right)\quad\forall\,(x,y,z)\in\mathbb{R}^{3}\,,
\end{equation}
and
\begin{equation}\label{M}
M(p,q):=\int_{[-T,T]\times[-\pi,\pi]}Q(X_{a}+p\mathbf{e}_{1}+q\mathbf{e}_{2})\cdot \partial_{t}X_{a}\wedge\partial_{\theta}X_{a}\,dt\,d\theta\quad\forall\,(p,q)\in\mathbb{R}^{2}\,.
\end{equation}

\begin{itemize}
\item[$(i)$]
If there exist $\bar \e >0$, a $C^0$-mapping $\varepsilon\mapsto(p_{\varepsilon},q_{\varepsilon})\in\mathbb{R}^{2}$ and a $C^1$-mapping $\varepsilon\mapsto\varphi_{\varepsilon}\in\mathcal{X}_{T}$ with $|\varepsilon|<\bar\varepsilon$, such that $\left.(p_\e, q_\e)\right|_{\e = 0} = (p_0, q_0)$, $\varphi_{0}=0$, and for every $\varepsilon\in(-\bar\varepsilon,\bar\varepsilon)$ the triple $(p_\e, q_\e, \f_\e)$ solves \eqref{M=H}, then $(p_{0},q_{0})$ is a critical point of $M$.

\item[$(ii)$]
If, in addition, $a\in\big(0,\frac12(\sqrt3-1)\big]$, $H_\e$ satisfies also $(H_2)$ and
\[
H_{\varepsilon}(x,y,z)=1+\varepsilon\widetilde{H}(x,y,z)\quad\forall\,(x,y,z)\in\mathbb{R}^{3}
\leqno{(H_{3})}
\]
with $\widetilde H \in C^2(\R^3)$, 
and $(p_{0},q_{0})\in\mathbb{R}^{2}$ is a nondegenerate critical point of $M$, then there exist $\bar \e >0$,  $C^1$-mappings $\varepsilon\mapsto(p_{\varepsilon},q_{\varepsilon})\in\mathbb{R}^{2}$, $\varepsilon\mapsto\varphi_{\varepsilon}\in\mathcal{X}_{T}$ with $|\varepsilon|<\bar\varepsilon$, such that $\varphi_{0}=0$, and for every $\varepsilon\in(-\bar\varepsilon,\bar\varepsilon)$ the triple $(p_\e, q_\e, \varphi_{\varepsilon})$ solves \eqref{M=H}.
\end{itemize}
\end{Theorem}
We point out that the mapping $M$ in \eqref{M} has the geometrical meaning of an algebraic, weighted volume of the body enclosed between $\Sigma_{a,T,p,q}$ and the planes of equation $z=\pm z_{a}(T)$, with a weight given by the perturbative part of $H_{\varepsilon}$, namely $\widetilde{H}(x,y,z)$. 

The proof of Theorem \ref{T:existence2} is based on the Lyapunov-Schmidt reduction, a technique which allows us to manage situations in which the linearized problem is defined by a Fredholm operator having a non null kernel. This is our case, since, when $a\in\big(0,\frac12(\sqrt3-1)\big]$ and $T\in\mathcal{T}_{a,1}\setminus\mathcal{T}_{a,0}$, the Jacobi operator $\mathfrak{L}_{a}$ turns out to have a two-dimensional kernel in $\mathcal{X}_{T}$.
By the Lyapunov-Schmidt method we can solve \eqref{M=H} for small $|\varepsilon|$ and for every $(p,q)\in\mathbb{R}^{2}$, apart from a couple of Lagrange multipliers, say $\lambda_{1,\varepsilon}(p,q)$ and $\lambda_{2,\varepsilon}(p,q)$, which correspond exactly to a basis of $\mathrm{ker}(\mathfrak{L}_{a})$ in the space $\mathcal{X}_{T}$. 

In order to remove them, we exploit the variational nature of \eqref{M=H}.
In fact, one can recognize that \eqref{M=H} is the Euler--Lagrange equation of a suitable energy functional in the space of parametric surfaces of annular type which are normal graphs of $\Sigma_{a,T,p,q}$ with fixed boundary, up to horizontal translations. 

Roughly speaking, $\lambda_{1,\varepsilon}(p,q)$ and $\lambda_{2,\varepsilon}(p,q)$ can be eliminated by taking variations of the energy along translations of $\Sigma_{a,T,p,q}$ in the directions $\mathbf{e}_{1}$ and $\mathbf{e}_{2}$, respectively. The gradient of the function $M$ comes out as leading term, with respect to $\varepsilon$, of such variations. For this reason, solutions to \eqref{M=H} can exist just in correspondence of critical points of $M$. Viceversa, nondegeneracy of critical points of $M$ is a sufficient condition for existence of solutions of \eqref{M=H}. 

We conclude with a couple of remarks: our results have some connections with those discussed in \cite{KoisoPalmerPiccione15} but with some differences: also in \cite{KoisoPalmerPiccione15}, degeneracy/nondegeneracy of compact sections of symmetric nodoids is studied, but in the context of bifurcation phenomena for surfaces with constant mean curvature $H_{0}$. There, the value $H_{0}$ is not fixed and can vary with $T$. Here, we fix $H_{0}=1$ and add a possibly variable perturbation on the prescribed mean curvature, and this makes the difference and arises new, distinct results. 

Finally, we observe that this study is part of a project which aims at constructing complete surfaces with prescribed, almost constant curvature, which cannot be detected in case of costant mean curvature. This kind of problem was first addressed in \cite{CaldiroliMusso}, looking for embedded compact surfaces with genus 1. However some mistakes were detected (see \cite{CIM, CaldiroliMusso2}) and a correct proof of the result stated in \cite{CaldiroliMusso} is not available, yet. The surfaces constructed in this work should be helping in this direction.

The paper is organized as follows: in Section \ref{S:preliminaries} we introduce some notation and we study the Jacobi operator $\mathfrak{L}_{a}$. Each of the next sections contains the proof of each theorem stated before, in the order. Finally, two appendices supplement this work: the first one about the Jacobi elliptic functions, the second one concerning the volume functional.

\section{Preliminaries}\label{S:preliminaries}
We start by listing some notations used throughout the paper. 

\begin{itemize}
\item We denote by $\{\mathbf{e}_1, \mathbf{e}_2, \mathbf{e}_3\}$ the canonical basis and by $\wedge$ the exterior product in $\R^3$.

\item If $(p, q) \in \R^2$ and $R>0$, then $B_R(p,q)$ is the ball of radius $R$ around $(p,q)$. If $(p,q) = (0, 0)$ we simply write $B_R$.

\item Given $T >0$, we define $\RT_T := [-T, T] \times [-\pi, \pi]$.

\item By $C$ we denote generic constants, whose value may change from line to line. 

\end{itemize}

\subsection{Surfaces of annular type and mean curvature}

By a surface of annular type we mean a surface $\Sigma$ in $\mathbb{R}^{3}$ parametrized by a sufficiently regular mapping $X\colon I\times\mathbb{S}^{1}\to\mathbb{R}^{3}$  where $I$ is an interval in $\mathbb{R}$. We identify $\mathbb{S}^{1}$ with $\mathbb{R}/2\pi\mathbb{Z}$ and we denote $(t,\theta)$ the pair of parameters. We always assume nice regularity, at least $C^{2}$, and, writing $X_{t}$ and $X_{\theta}$ the derivatives of $X$ with respect to $t$ and $\theta$, respectively, we also assume that $X_{t}\wedge X_{\theta}\ne 0$ at every $(t,\theta)\in I\times[-\pi,\pi]$. The normal vector to $\Sigma$ at a given point $X\in\Sigma$ is 
\[
{N}:=\frac{X_{t}\wedge X_{\theta}}{\left|X_{t}\wedge X_{\theta}\right|}~\!.
\]
The mean curvature of $\Sigma$ at $X\in\Sigma$, denoted by $\mathfrak{M}(X)$, is expressed in terms of the coefficients of the first and second fundamental form, in Gaussian notation,
\[
\begin{array}{c}
{\E}=|X_{t}|^{2}\,,\quad {\F}=X_{t}\cdot X_{\theta}\,,\quad {\G}=|X_{\theta}|^{2}\,,\\
{\L}=X_{tt}\cdot {\normal}\,,\quad {\M}=X_{t\theta}\cdot {\normal}\,,\quad {\mathcal{N}}=X_{\theta\theta}\cdot {\normal}\,,
\end{array}
\]
as
\[
\mathfrak{M}=\frac{{\E}{\mathcal{N}}-2{\F}{\M}+{\G}{\L}}{2({\E}{\G}-{\F}^{2})}\,.
\]

A surface of revolution $\Sigma$ in $\mathbb{R}^{3}$ with axis of revolution given by the $z$-axis is the image of a mapping $X\colon I\times \mathbb{S}^1\to\mathbb{R}^3$ of the form
\[
X(t,\theta)=\left[\begin{array}{c}x(t)\cos\theta\\ x(t)\sin\theta\\ z(t)\end{array}\right]\quad\forall (t,\theta)\in I\times[-\pi,\pi]\,,
\]
where $x\colon I\to (0,\infty)$ and $z\colon I\to\mathbb{R}$ are mappings of class $C^{2}$. The curve $t\mapsto(x(t),0,z(t))$, with $t\in I$, defines the generatrix of the surface $\Sigma$. 


\subsection{The Jacobi operator}

Fixing $a \in (-1,0) \cup(0, \infty)$, the Jacobi operator $\mathfrak{L}_{a}$ corresponding to the Delaunay surface $\Sigma_{a}$ is defined as the linearized mean curvature operator around the parameterization \eqref{Xa} with respect to normal variations. More precisely, for every $\varphi\in C^{2}(\mathbb{R}\times\mathbb{S}^{1})$ we set
\begin{equation}\label{La-def}
\mathfrak{L}_{a}\varphi:=\left.\frac{\partial}{\partial s}[\mathfrak{M}(X_{a}+s\varphi N_{a})]\right|_{s=0}\quad\text{where}\quad N_{a}:=\frac{(X_{a})_{t}\wedge(X_{a})_{\theta}}{|(X_{a})_{t}\wedge(X_{a})_{\theta}|}\,.
\end{equation}
We point out that the definition of $\mathfrak{L}_{a}\varphi$ is well posed. Indeed, 
for every $\varphi\in C^{2}(\mathbb{R}\times\mathbb{S}^{1})$, we have 
\begin{equation}\label{eq:extcond}
(X_{a}+\varphi N_{a})_{t}\wedge(X_{a}+\varphi N_{a})_{\theta}\cdot N_{a}=x_{a}^{2}\left[1-2\varphi+\frac{\varphi^{2}}{x_{a}^{2}}\left(x_{a}^{2}-\frac{\gamma_{a}^{2}}{x_{a}^{2}}\right)\right]\,.
\end{equation}
Hence $[(X_{a}+s\varphi N_{a})_{t}\wedge(X_{a}+s\varphi N_{a})_{\theta}](t,\theta)\ne 0$ for $|s|$ small enough, depending on $a$, $\varphi$, $(t,\theta)$, which allows us to write $\mathfrak{M}(X_{a}+s\varphi N_{a})(t,\theta)$ for $|s|$ small, and to compute its derivative at $s=0$. 

\begin{Proposition}\label{L:Jacobi}
It holds 
\[
\mathfrak{L}_{a}\varphi=\frac{1}{2x_{a}^{2}}\left(\Delta\varphi+2p_{a}\varphi\right)\,,
\]
where
\[
p_{a}(t):=x_{a}^{2}(t)+\frac{\gamma_{a}^{2}}{x_{a}^{2}(t)}\,,\qquad\gamma_{a}:=a(a+1)\,,
\] 
and $\Delta=\partial_{tt}+\partial_{\theta\theta}$.
\end{Proposition}
We refer to \cite[Lemma 3.1]{CaldiroliMusso} for a proof. 
\begin{Remark}\label{R:paprop}
By definition, $p_a$ is even and $2\tau_a$-periodic. Moreover, by means of \eqref{eq:xaza}, \eqref{eq:isothermal} we see that $p_a$ satisfies the equation
\[
p_a'' = 2\left[2(1+2\gamma_{a})p_a - 3p_a^2 + 4\gamma_{a}^2\right]\,.
\]
Then, recalling that $\min\{1+a, |a|\} \leq x_a(t) \leq \max\{1+a, |a|\}$ for every $t \in \R$, where $x_a( 0) = 1+a$ and $x_a(\tau_a) = |a|$, by a standard argument on autonomous second order ODEs we infer that the fundamental period of $p_a$ is actually $\tau_a$ and that
\begin{equation}\label{eq:paprop}
\max_{t\in \R}p_a = p_a(0 ) = 1+2\gamma_{a}\,, \qquad \min_{t \in \R}p_a = p_a\left(\frac{\tau_a}{2}\right) =  2|\gamma_{a}|\,.
\end{equation}
\end{Remark}

\subsection{Fundamental solutions of the linearized problem}
Fix $a \in (-1, 0)\cup (0,\infty)$ and $T>0$. Let $\f \in C^2([-T, T]\times \S^1)$ be such that
\[
\begin{gathered}
\f (t, \theta) = w(t) \cos(j\theta)\, \text{ for some } j \in \N\,,\\
\f(\pm T, \cdot ) = 0, \quad \f(t, \cdot) = \f(-t, \cdot)\,.
\end{gathered}
\]
Then $\mathfrak L_a\f = 0$ if and only if $w \in C^{2}([-T, T])$ is even and solves 
\begin{equation}\label{eq:j-eq}
\begin{cases}
w'' = (j^2-2p_a)w & \text{ in }\ (-T, T)\\
w(\pm T) = 0\,.
\end{cases}
\end{equation}

Our goal is to determine whether, for any fixed couple $a, j$,  there exist some $T>0$ such that \eqref{eq:j-eq} admits a nontrivial even solution. To this aim, it seems convenient to begin by studying the nodal set of the solutions to the more general equation
\begin{equation}\label{eq:geneq}
u'' = q u\quad \text{ in }\ \R\,,
\end{equation}
without imposing any boundary condition. Here, we assume that $q \in C^0(\R)$ is even and $\tau$-periodic. By well known results every solution to \eqref{eq:geneq} is a linear combination of the principal fundamental solutions $w, v \in C^2(\R)$ which satisfy
\[
\begin{cases}
w(0) = 1\\
w'(0) = 0\,,
\end{cases}
\qquad 
\begin{cases}
v(0) = 0\\
v'(0) = 1\,.
\end{cases}
\]
Since $q$ is even, we have that $w$ is even and $v$ is odd.

\begin{Remark}\label{R:Floquet}  
Applying Floquet Theory to Hill's equation (see e.g. \cite[Section 2.4.2]{Chicone}) one obtains a partial characterization of $w$, $v$. Indeed, there exist $r, s:\R \to \R$, with $r$ even and $s$ odd, such that exactly one of the following alternatives holds true:
\begin{itemize}
\item[$(i)$] $r,s$ are both either $\tau$-periodic or $2\tau$-periodic, and there exists $\mu >0$ such that
 \[
\begin{aligned} 
&w (t)= r(t)\cosh(\mu t) + s(t)\sinh(\mu t)\\
&v (t)= r(t)\sinh(\mu t) + s(t)\cosh(\mu t)\,;
\end{aligned}
\]

\item[$(ii)$]
$r,s$ are both either $\tau$-periodic or $2\tau$-periodic, and there exists $\lambda \in \R$ such that
 \[
\begin{aligned} 
&w(t)= r(t)\\
&v(t) = s(t) + \lambda t r(t)\,;
\end{aligned}
\]

\item[$(iii)$]
$r,s$ are both either $\tau$-periodic or $2\tau$-periodic, and there exists $\lambda \in \R$ such that
 \[
\begin{aligned} 
&w(t) = r(t) + \lambda t s(t) \\
&v(t) = s(t)\,;
\end{aligned}
\]
\item[$(iv)$]
$r,s$ are both $\tau$-periodic and there exists $0<\sigma<\frac{\pi}{\tau} $ such that
\[
\begin{gathered}
w(t)  = r(t) \cos(\sigma t) - s(t) \sin(\sigma t)\\
v(t) = r(t)\sin(\sigma t) + s(t)\cos(\sigma t)\,.
\end{gathered}
\]
\end{itemize}

\end{Remark}

\begin{Lemma}\label{L:genq}
Let $q \in C^0(\R)$ be a even $\tau$-periodic function, and let $w$ be the even principal fundamental solutions of $u''= qu$ in $\R$. 
\begin{enumerate}
\item[$(i)$] If
\begin{equation}\label{eq:meancond}
 \frac{1}{\tau}\int_{0}^{\tau} q \, dt\ < 0\,,
\end{equation}
then there exists $T>0$ such that $w(\pm T) = 0$.
\item[$(ii)$] If $w(\alpha) = w(\beta) = 0$ for some $\alpha,\beta \in \R$ with $\alpha<\beta$, then
\[
q^- \not \equiv 0\,,
\]
where $q^- = \max\{-q, 0\}$ is the negative part of $q$. Moreover
\[
\beta-\alpha \geq \pi\left(\underset{t \in [0, \tau]}{\max} q^-(t)\right)^{-\frac{1}{2}}\,.
\]
\end{enumerate}
\end{Lemma}
\begin{Proof} Let us prove $(i)$. Assume that \eqref{eq:meancond} holds, and suppose by contradiction that $w \neq 0$ for every $t \in \R$, that is, $w >0$ on $\R$ (because $w(0)= 1$). 
Fix $s >0$. Since $w>0$ and $w'(0) = 0$, from \eqref{eq:geneq} we find
\[
\int_0^{s} q(\xi) \, d \xi\, = \int_0^{s} w^{-1}(\xi)w''(\xi)\,d\xi = w(s)^{-1}w'(s)+ \int_0^{s} \left(w^{-1}(\xi)w'(\xi)\right)^2\, d\xi\,.
\] 
Let $k_s \in \N$ be the integer part of $\tau^{-1}s$, so that $s \in [k_s\tau, (k_s+1)\tau)$. Being $q$ $\tau$-periodic, we get
\[
w(s)^{-1}w'(s) < \sum_{i=0}^{k_s-1} \int_{i\tau}^{(i+1)\tau} q(\xi)\, d\xi + \int_{k_s\tau}^s q(\xi)\, d\xi \leq \tau I k_s + (s - k_s \tau)\max_{\xi \in [0,\tau]}q(\xi)\,,
\]
where we set $I := \tau^{-1}\int_{0}^{\tau} q(\xi)\, d\xi$.
Thus, using that $s < (k_s+1)\tau$, we obtain 
\begin{equation}\label{eq:logprop}
w(s)^{-1}w'(s) <\tau(I k_s + \max_{\xi \in [0,\tau]}q(\xi))\,, \quad \text{ for any } s>0.
\end{equation}

Let $ t >0$. Integrating \eqref{eq:logprop} over $[0, t]$ we find
\[
0 <w( t)  \leq e^{\tau\left(I \int_0^{t} k_s ds  + t\underset{\xi \in [0, \tau]}{\max}q(\xi)\right)} \quad  \text{ for any } t > 0\,,
\]
where we used that $w(0) = 1$. Now, an explicit computation shows that
\[
I \int_0^{t} k_s ds  + t  \max_{\xi \in [0,\tau]}q(\xi) = \frac{I}{2}(\tau k_{t})^2(1 + o(1)) \quad \text{ as }t \to \infty, 
\]
thus, since by assumption $I <0$, we infer that 
\begin{equation}\label{eq:asym}
w(t) \to 0 \quad \text{ as }\quad t \to \infty. 
\end{equation}

To complete the proof of $(i)$ we argue by exhaustion. Consider the sequence $a_n :=w(2\tau n)$, $n \in \N$ and suppose that $w$ is as in Remark \ref{R:Floquet}, $(i)$. Since $r, s$ are at least $2\tau$-periodic, $s$ is odd, and $w(0) = 1$, then $a_n = \cosh(2\tau\mu n)$, which implies that $\limsup_{t \to \infty}w = \infty$, contradicting \eqref{eq:asym}. Similarly, if $w$ is as in Remark \ref{R:Floquet}, $(ii)-(iii)$, we get that $a_n = 1$, again a contradiction. Finally, let $w$ be as in Remark \ref{R:Floquet}, $(iv)$. Arguing as before, we have 
\[
a_n = r(2\tau n) \cos(2\tau \sigma n) + s(2\tau n) \sin(2\tau \sigma n) =  \cos(2\tau \sigma \,n)\,. 
\]
We now distinguish two cases: if $\frac{\tau\sigma}{\pi} \in \mathbb{Q}$, then we can extract a subsequence such that $a_n  = 1$, again a contradiction. On the other hand if  $\frac{\tau\sigma}{\pi} \in \R\setminus \mathbb{Q}$ then $a_n$ is dense in $[-1,1]$. In particular there exists some $n_0$ such that $a_{n_0} = w(2\tau n_0) <0$.
Hence, also in this case we obtain a contradiction, since we are assuming $w >0$. The proof of $(i)$ is complete. 

\smallskip

Let now prove $(ii)$. Testing \eqref{eq:geneq} with $w$ and integrating by parts we obtain
\[
0 = \int_\alpha^\beta |w'|^2 dt + \int_\alpha^\beta q|w|^2 dt\,.
\]
By the one dimensional Poincar\'e inequality we have
\[
\int_\alpha^\beta |w'|^2 dt \geq \left((\beta-\alpha)^{-1}\pi\right)^2\int_\alpha^\beta |w|^2 dt
\]
thus we find 
\[
0 \geq \left[ \left((\beta-\alpha)^{-1}\pi\right)^2 + \min_{[0,\tau]}q \right]\int_\alpha^\beta |w|^2 dt\,,
\]
that is,
\[
0< \left((\beta-\alpha)^{-1}\pi\right)^2 \leq - \min_{[0,\tau]}q  = \max_{[0,\tau]} q^- \,,
\]
and $(ii)$ easily follows. The Lemma is proved.  \qed
\end{Proof}

We now focus on  the equation (compare with \eqref{eq:j-eq})
\begin{equation}\label{eq:ajeq}
u'' = (j^2-2p_a)u \quad \text{ in }\ \R\,.
\end{equation}
As next Lemma shows, solutions to \eqref{eq:ajeq} are completely known when $j=0,1$, see \cite[Lemma 4.3]{CaldiroliMusso}.
\begin{Lemma} \label{P:ker}
Let $a\in(-1,0)\cup(0,\infty)$.
\begin{enumerate}
\item[$(i)$] Set
\[
w_{a,0} :=
\begin{cases}\displaystyle
 -\frac{z'_{a}}{x_{a}}\frac{\partial x_{a}}{\partial a}+\frac{x'_{a}}{x_{a}}\frac{\partial z_{a}}{\partial a} & \text{if } a \neq -\frac12\,,\\
-\cos(t) &\text{if } a = -\frac12\,,
\end{cases}
\qquad\qquad v_{a,0} := 
\begin{cases}\displaystyle
\frac{x'_{a}}{x_{a}}\,, &\text{if }  a \neq -\frac{1}{2}\,,\\
-\sin(t)&\text{if } a = -\frac12\,.
\end{cases}
\]
Then every solution to $u'' = -2p_a u$ in $\R$ is a linear combination of $w_{a,0}, v_{a,0}$. Moreover, $w_{a,0}$ is even and such that $w_{a,0}(0) = -1$ while $v_{a,0}$ is odd. 
\item[$(ii)$] Set
\[
{w_{a,1} := \frac{z'_a}{x_a}}\,, \qquad v_{a,1} := \frac{x_ax_a'+z_az_a'}{x_a}\,.
\]
Then every solution to $u'' = (1 -2p_a )u$ in $\R$ is a linear combinations of $w_{a,1}, v_{a,1}$. Moreover, $w_{a,1}$ is even and such that $w_{a,1}(0) = 1$ while $v_{a,1}$ is odd. 
\end{enumerate}
\end{Lemma}

In general, given $j \in \N$ we denote by $w_{a,j}$ the unique solution to the Cauchy problem
\begin{equation}\label{eq:wv}
\begin{cases}
w_{a,j}'' = (j^2-2p_a)w_{a,j} & \text{ in }\ \R\\
w_{a,j}(0) = (-1)^{j+1}\\
w_{a,j}'(0)= 0 
\end{cases}
\end{equation}
Clearly, since $j^2-2p_a$ is even, $w_{a,j}$ is even. Moreover, any even solution to \eqref{eq:ajeq} is proportional to $w_{a,j}$. 
As a straightforward consequence of Lemma \ref{L:genq} we obtain the following. 
\begin{Corollary}\label{L:zeridist}
Let $a \in (-1, 0) \cup (0, \infty)$ and  $j\geq 2$. 
\begin{enumerate}
\item[$(i)$] If 
\begin{equation}\label{eq:nochangecond}
j^2 \geq  2(1+2\gamma_a) \,,
\end{equation}
then $w_{a,j}(t) \neq 0$ for every $t \in \R$.

\item[$(ii)$] If 
\[
j^2 < 4|\gamma_a| \,,
\]
then there exists $T >0$ such that $w_{a,j}(\pm T)= 0$.
\end{enumerate}
\end{Corollary}

\begin{Proof} Let us prove $(i)$. Assume by contradiction that there exists $T >0$ such that $w_{a,j}(T) = 0$. Since $w_{a,j}$ is even, we have $w_{a,j}(T) = w_{a,j}(-T) = 0$. Hence, by Lemma \ref{L:genq}, $(ii)$, with $q = j^2-2p_a$ we infer that $j^2 - 2p_a(t_*) < 0$ for some $t_* \in [0, \tau_a]$.  On the other hand, by \eqref{eq:paprop} and \eqref{eq:nochangecond} it holds
\[
j^2-2p_a \geq j^2 - 2(1+2\gamma_a) \geq 0\,,
\]
a contradiction. 

\smallskip
 
As for $(ii)$, again by Remark \ref{R:paprop} we infer
\[
\frac{1}{\tau_a}\int_{0}^{\tau_a} j^2-2p_a d t \leq j^2 - 4|\gamma_a|\,.
\]
Applying Lemma \ref{L:genq}, $(i)$ with $q = j^2-2p_a$ we get the result.  \qed
\end{Proof}

Next Lemma provides an explicit formula for $w_{a,0}$.
\begin{Lemma}\label{prop:explicitformwa0+}
Let $a\in\left(-1,0\right) \cup(0,\infty) \setminus \left\{-\frac12\right\}$. For any $t\in\R$ it holds
\begin{equation}\label{eq:formula}
w_{a,0}(t)=\frac{1}{1+2a}\left[1+2\gamma_a -2x^2_a(t)+\frac{x_a'(t)}{x_a(t)} (2z_a(t) - t )\right]\,. 
\end{equation}
\end{Lemma}
\begin{Proof} The result can be checked directly, by showing that both $w_{a,0}$ and the right-hand side of \eqref{eq:formula} solve $u'' = -2p_a u $ in $\R$ with the initial conditions $u(0) = -1$, $u'(0) = 0$. 

Nevertheless, we also provide a constructive proof. To this goal, it is useful to recall the relation between $x_a$ and Jacobi elliptic functions. By using and adapting \cite[Proof of Lemma 2.3]{CIM}, one can see that for any $t \in \R$ it holds 
\begin{equation}\label{eq:formulaxa}
x_a(t)=
\begin{cases}
\displaystyle (1+a)\, \mathrm{dn}((1+a)t| m_a)\,, & \text{ if } a \in \left(-\frac{1}{2}, 0\right) \cup(0, \infty)\,,\\
\displaystyle |a|\, \dn(|a|(t-\tau_a)| m_a)\,, & \text{ if } a \in \left(-1, -\frac{1}{2}\right)\,,\\\end{cases}
\end{equation}
where $\dn(s| m)$ is the Delta Amplitude function (see Appendix \ref{A:Jacobi}) with {parameter}
\[
m_a:=
\begin{cases}
\displaystyle 1-\frac{a^2}{(1+a)^2}  & \text{ if } a \in \left(-\frac{1}{2}, 0\right) \cup(0, \infty)\,,\\[10pt]
\displaystyle 1-\frac{(1+a)^2}{a^2}& \text{ if } a \in \left(-1, -\frac{1}{2}\right)\,.
\end{cases}
\]
Moreover, we recall that 
\[
K(m) = \int_0^{\pi/2}\frac{1}{\sqrt{1-m\sin^2 \theta}}\,d\theta \quad \text{  }\quad E(m) =\int_{0}^{\pi/2}\sqrt{1-m\sin^2\theta}\, d\theta\,, \qquad m \in [0,1)\,,
\]
are, respectively, the complete elliptic integrals of the first and second kind and that it holds (again, one can argue as in \cite[Proof of Lemma 2.3]{CIM})
\begin{equation}\label{eq:tauK}
\tau_a =
\begin{cases}
\displaystyle \frac{K(m_a)}{1+a}& \text{ if } a \in \left(-\frac{1}{2}, 0\right) \cup(0, \infty)\,,\\[10pt]
\displaystyle \frac{K(m_a)}{|a|} & \text{ if } a \in \left(-1, -\frac{1}{2}\right)\,.
\end{cases}
\end{equation}
Finally, using \eqref{eq:formulaxa}, \eqref{eq:tauK} and \eqref{eq:dnintcei}, we easily see that 
\begin{equation}\label{R:ha}
\int_0^{\tau_a} x_a^2 ds = 
\begin{cases}
\displaystyle (1+a) E(m_a) & \text{ if } a \in \left(-\frac{1}{2}, 0\right) \cup(0, \infty)\,,\\
\displaystyle |a|E(m_a)\,,& \text{ if } a \in \left(-1,-\frac{1}{2}\right) \,.
\end{cases}
\end{equation}

We are now in the position to compute the derivative of $x_a$ with respect to $a$.  First we assume that $a\in\left(-\frac{1}{2},0\right) \cup(0,\infty)$, and use the first representation in \eqref{eq:formulaxa}. Thus we get
\begin{equation}\label{eq:dxdapass1}
\begin{aligned}
\frac{\partial x_a}{\partial a}(t)= \mathrm{dn}((1+a)t|m_a) &+ (1+a)\,t\, \left.\frac{\partial}{\partial s}\left[ \mathrm{dn}(s|m)\right]\right|_{(s,m)=((1+a)t, m_a)} \\
&- \frac{a}{1+2a}\,2m_a \left.\frac{\partial}{\partial m}\left[ \mathrm{dn}(s|m)\right]\right|_{(s,m)=((1+a)t, m_a)}\,.
\end{aligned}
\end{equation}
From \eqref{eq:formulaxa} we have
\begin{equation}\label{eq:dertxa}
\left.\frac{\partial}{\partial s}\left[ \mathrm{dn}(s|m)\right]\right|_{(s,m)=((1+a)t, m_a)} = \frac{x_a^\prime(t)}{(1+a)^2} \,.
\end{equation}
Moreover, using \eqref{eq:dndk}, \eqref{eq:dertxa} and recalling that $1-m_a = (1+a)^{-2}a^2$, we find
\begin{equation}\label{eq:dndkint}
\begin{aligned}
2m_a\left.\frac{\partial}{\partial m}\left[ \mathrm{dn}(s|m)\right]\right|_{(s,m)=((1+a)t, m_a)} &= \frac{x_a(t)}{1+a} - \frac{1+a}{x_a(t)} + \frac{x_a'(t)}{1+a}t\\
& - \frac{\left(x_a'(t)\right)^2}{a^2(1+a)x_a(t)} - \frac{x_a'(t)}{a^2}\int_0^{(1+a)t}\mathrm{dn}(\tau |m_a)^2 d\tau \,.
\end{aligned}
\end{equation}
Now, from \eqref{eq:xaza} we have that
\begin{equation}\label{eq:expressza}
z_a(t)=\int_0^tx_a^2(s)\,ds -\gamma_a t\,.
\end{equation}
As a consequence, using also the change of variables $\tau = (1+a)s$ we get
\begin{equation}\label{eq:dxdapass2}
\int_0^{(1+a)t}\mathrm{dn}(\tau|m_a)^2 d\tau = (1+a)\int_0^t (\mathrm{dn}((1+a)s|m_a)^2 ds = \frac{1}{1+a}\int_0^tx_a^2(s)\,ds = \frac{z_a(t) + \gamma_a t}{1+a}\,.
\end{equation}

Putting together \eqref{eq:dxdapass1}--\eqref{eq:dndkint} and \eqref{eq:dxdapass2} we find
\[
\begin{aligned}
\frac{\partial x_a}{\partial a}(t)&= \frac{x_a(t)}{1+a} +\frac{x_a'(t)}{1+a}t\\ 
&- \frac{ a}{1+2a}\,\left( \frac{x_a(t)}{1+a} - \frac{1+a}{x_a(t)} + \frac{x_a'(t)}{1+a}t - \frac{\left(x_a'(t)\right)^2}{a^2(1+a)x_a(t)}
 - \frac{x_a'(t)(z_a(t) + \gamma_a t)}{ a^2(1+a)}\right)\,.
\end{aligned}
\]
Exploiting \eqref{eq:xaza} and \eqref{eq:isothermal} we have the chain of equalities
\begin{equation}\label{eq:chain}
(x_a')^2 = x_a^2 - (z_a')^2 = (1+2\gamma_a)x_a^2 - x_a^4 - \gamma_a^2 =(1+\gamma_a-z_a')x_a^2 - \gamma_a^2 \,,
\end{equation}
Thus, using \eqref{eq:chain} and performing some standard computations we get 
\begin{equation}\label{eq:dxda}
\frac{\partial x_a}{\partial a}(t)=  \frac{1}{\gamma_a(1+2a)}\left[(1+2\gamma_a) x_a(t)+ 2\gamma_ax_a'(t)\,t+ x_a'(t)z_a(t) - x_a(t)z'_a(t)\right]\,,
\end{equation}
which conclude the computation of $\frac{\partial{x_a}}{\partial a}$ when $a \in \left(-\frac{1}{2},0\right)\cup(0,\infty)$. 

Let now $a \in \left(-1, -\frac12\right)$. Here we have to use the second representation in \eqref{eq:formulaxa}.
Nevertheless, we can argue as before, with some obvious modifications and paying attention to the further technical issue due to the presence of $\tau_a$ in the representation of $x_a$. We omit the details: we limit ourselves to point out that, in order to overcome the additional difficulty, one needs to  take \eqref{R:ha} into account and use that 
\[
\frac{\partial \tau_a}{\partial a} = \frac{(1+a)K(m_a) + aE(m_a)}{(1+2a)\gamma_a}\,,
\]
which can be obtained by direct computation by means of \eqref{eq:tauK} and, for instance, \cite[Eq. 19.4.1]{DLMF}. In the end, also in this case \eqref{eq:dxda} holds. 
\medskip

Next we compute the derivative of $z_a$ with respect to $a$. To ease the notation, from now on we avoid to explicitly state the variables when denoting a function. Differentiating \eqref{eq:expressza} with respect to $a$, exploiting \eqref{eq:chain}, \eqref{eq:dxda} and integrating by parts we get
\begin{equation}\label{eq:derzarispettoa}
\begin{aligned}
\frac{\partial z_a}{\partial a}&=2\int_0^tx_a\frac{\partial x_a}{\partial a}\,ds -(1+2a) t\\
&= \frac{2+5\gamma_a}{\gamma_a(1+2a)}\int_0^t x_a^2\,ds + \frac{1}{\gamma_a(1+2a)}\left(2\gamma_a x_a^2\,t+x_a^2z_a - 3\int_0^tx^4_a\,ds\right)  -(1+2a) t\,.
\end{aligned}
\end{equation}

To compute the integral of $x_a^4$ we first notice that by \eqref{eq:xaza} it holds
\[
\int_0^{t}x_ax''_a ds = (1+2\gamma_a )\int_0^{t}x^2_a ds - 2 \int_0^{t}x^4_a ds\,.
\]
On the other hand, integrating by parts and using \eqref{eq:chain} we get
\[
\int_0^{t}x_ax''_a ds = x_ax_a' - \int_0^{t}(x'_a)^2 ds = x_ax_a' -(1+2\gamma_a) \int_0^{t}x_a^2 ds +\int_0^t x_a^4\,ds + \gamma_a^2\,t\,.
\]
Comparing these expressions we find 
\begin{equation}\label{eq:x4int}
3\int_0^t x_a^4\,ds  =2 (1+2\gamma_a )\int_0^{t}x^2_a ds - x_ax_a' - \gamma_a^2\,t\,.
\end{equation}

Plugging \eqref{eq:x4int} into \eqref{eq:derzarispettoa} we find
\[
\frac{\partial z_a}{\partial a}= \frac{1}{1+2a}\int_0^t x_a^2\,ds + \frac{1}{\gamma_a(1+2a)}\left(2\gamma_a x_a^2\,t+x_a^2z_a +x_ax_a' + \gamma_a^2\,t\right)  -(1+2a) t\,.
\]
Finally, using \eqref{eq:expressza} and \eqref{eq:chain}, after some elementary computations and we get
\begin{equation}\label{eq:dzda}
\frac{\partial z_a}{\partial a} =\frac{1}{\gamma_a(1+2a)}\left[2\gamma_az_a + 2\gamma_a z_a'\,t- \gamma_a t +x_ax_a'+ z_az'_a \right] \,.
\end{equation}

We are in position to complete the proof. Putting \eqref{eq:dxda} and \eqref{eq:dzda} in the expression of $w_{a,0}$ (see Lemma \ref{P:ker}) and regrouping the terms we infer that
\[
w_{a,0}=\frac{1}{1+2a}\left[\frac{(z'_a)^2+(x_a')^2-(1+2\gamma_a) z'_a}{\gamma_a} +\frac{x_a'}{x_a} (2z_a - t )\right]\,.
\]
The conclusion easily follows by \eqref{eq:xaza} and \eqref{eq:isothermal}.\qed
\end{Proof}

\section{Proof of Theorem \ref{T:Ta}}

In this Section we prove Theorem \ref{T:Ta}. We begin with a preliminary result. 
Recalling the definition of $w_{a,j}$ (see Lemma \ref{P:ker} and \eqref{eq:wv}) we can rewrite the sets \eqref{Taj} as 
\[
\T_{a,j}=\{t>0 \mid w_{a,j}(t)=0\}\,, \quad j \in \N\,.
\]

Concerning $\T_{a,0}$, $\T_{a,1}$ and their relationship, we have the following.
\begin{Lemma}\label{L:nod0}
Let $a \in (-1, 0) \cup (0, \infty)$.
\begin{enumerate}
\item[$(i$)] If  $a = -\frac{1}{2}$, then $\T_{a,0} = \left\{\left(k+\frac{1}{2}\right) \pi\right\}_{k=0}^\infty$. If $a \neq -\frac{1}{2}$, then $\T_{a,0} = \{T_{a,k}\}_{k=0}^\infty$, where
\[
k\tau_a<T_{a,k} < \left(k+\frac{1}{2}\right)\tau_a, \quad \text{ for every }k \in \N\,.
\]

\item[$(ii)$] If $a \in (-1, 0)$, then $\T_{a,1} = \emptyset$. If $a \in (0, \infty)$, then $\T_{a,1} = \left\{ \left(k+ \frac{1}{2}\right)\tau_a \right\}_{k=0}^\infty$.

\item[$(iii)$] $\T_{a,0}\cap \T_{a,1} = \emptyset$. 
\end{enumerate}
\end{Lemma}
\begin{Proof} We begin by proving $(ii)$ and $(iii)$. Using the definition of $w_{a,1}$ (see Lemma \ref{P:ker}) and  \eqref{eq:xaza}, it is easy to see that $w_{a,1} > 0$ when $a \in (-1, 0)$. Therefore $\T_{a,1}=\emptyset$ and thus $\T_{a,0}\cap \T_{a,1}= \emptyset$ is trivial in this case.

Assume now that $a>0$. Using Remark \ref{R:paprop} one can see that both $w_{a,1}(t)$ and $ -w_{a,1}(t + \tau_a)$ solve
$u'' = (1-2p_a)u$ in $\R$ with initial conditions $u(0) = -1$, $u'(0) = 0$. Therefore $w_{a,1}(t) + w_{a,1}(t + \tau_a) = 0$ for every $t\in \R$. By taking $t = -\frac{\tau_a}{2}$ and using the symmetry of $w_{a,1}$ we get 
\[
w_{a,1} \left(\frac{\tau_a}{2}\right) = w_{a,1} \left(-\frac{\tau_a}{2}\right) =0\,,
\]
thus from the periodicity of $w_{a,1}$ we infer that  $\left\{ \left(k+ \frac{1}{2}\right)\tau_a \right\}_{k=0}^\infty \subseteq \T_{a,1}$. On the other hand, since
\[
w'_{a,1}(t) = x'_a(t)\left(1 + \frac{\gamma_a}{x_a^2(t)}\right)\,,
\]
we have that $w_{a,1}$ is monotone decreasing in $(0, \tau_a)$ and monotone increasing in $(\tau_a, 2\tau_a)$. Hence, since $w_{a,1}(0) = 1$, it follows that $w_{a,1}$ has at most one zero in $(0, \tau_a)$, and at most one zero in $(\tau_a, 2\tau_a)$. Now, observing that $w_{a,1}(\frac{\tau_a}{2}) = w_{a,1}(\frac{3\tau_a}{2}) = 0$, and using the periodicity of $w_{a,1}$ we infer that $\T_{a,1} \subseteq \left\{ \left(k+ \frac{1}{2}\right)\tau_a \right\}_{k=0}^\infty$ and $(ii)$ follows.

Let us set $t_k := \left(k+\frac{1}{2}\right)\tau_a$, $k\in \N$. Since $w_{a, 1}(t_k) =0$ and  $x_a(t) >0$ for any $t \in \R$, we immediately get that $z_a'(t_k) = 0$ and $x_a^2(t_k) = \gamma_a$. Thus, by \eqref{eq:isothermal} we find that $x_a(t_k) = |x_a'(t_k)|$. Recalling that $x_a'(t)<0$ whenever $t \in \bigcup_{n\in\N}(2n\tau_a, (2n+1)\tau_a)$, we infer that
\[
\frac{x_a'(t_k)}{x_a(t_k)} = (-1)^{k+1}\,.
\]
Then, by Lemma \ref{P:ker} and \eqref{eq:dzda} we obtain 
\[
w_{a,0}(t_k) = (-1)^{k+1}\frac{x_a(t_k)x'_a(t_k) +  2\gamma_a z_a(t_k) -\gamma_at_k}{\gamma_a(1+2a)} = \frac{(-1)^{k+1}}{1+2a}\left(\frac{x_a'(t_k)}{x_a(t_k)}+ 2z_a(t_k) - t_k\right)\,.
\]

Let us define 
\[
f(t) = \frac{x_a'(t)}{x_a(t)} +  2 z_a -t\,, \quad t \in \R\,.
\]
By the previous discussion, we have that $w_{a,0}(t_k) = 0$ if and only if $f(t_k) = 0$. 
Using \eqref{eq:xaza} and \eqref{eq:isothermal} we see that $f(0)= 0$ and  
\[
f'(t) = p_a(t) - (1+2\gamma_a)\,.
\]
Hence, thanks to Remark \ref{R:paprop} we get that $f'(t) \leq 0$ for every $t \in \R$, and $f'(t) = 0$ if and only if $t = n \tau_a$ for some $n \in\Z$. As a consequence we have that $f(t) =0$ if and only if $t=0$.  
Therefore $w_{a,0}(t_k)\neq 0$ for any $k \in \N$ and thus $\T_{a,0} \cap \T_{a,1} = \emptyset$. The proof of $(iii)$ is complete. 

\smallskip

As for $(i)$, if $a = -\frac{1}{2}$ the conclusion follows immediately by Lemma \ref{P:ker}. Let then $a\in (-1,0)\cup(0,+\infty)$ with $a\neq -\frac{1}{2}$. As a further consequence of the proof of $(ii)$ we have that
\[
w_{a,0} \left(\left(k+\frac{1}{2}\right)\tau_a\right) = (-1)^{k}\left|w_{a,0} \left(\left(k+\frac{1}{2}\right)\tau_a\right)\right|\,, \quad \text{ for any }k \in \N\,.
\]
Moreover, by means of Lemma \ref{prop:explicitformwa0+} we can compute
\[
w_{a,0}(k\tau_a)=\frac{1}{1+2a}\left[1+2\gamma_a -2x^2_a(k\tau_a)\right] = (-1)^{k+1}\,,\quad \text{ for any }k\in \N\,,
\]
where we used that $x(2n\tau_a) = 1+a$ and $x((2n+1)\tau_a) = |a|$ for any $n\in \Z$. Thus we infer that
\[
w_{a,0}(k\tau_a)w_{a,0} \left(\left(k+\frac{1}{2}\right)\tau_a\right) <0 \quad \text{ for any }k\in\N\,.
\]
As a consequence, there exists $T_{a,k} \in \left(k\tau_a, \left(k+\frac{1}{2}\right)\tau_a\right)$ such that $w_{a,0}(T_{a,k})=0$ for any $k\in\N$.

Finally, let us recall that $w_{a,0}$ and $v_{a,0}$ are two linearly independent solutions to \eqref{eq:ajeq} with $j=0$ (see Lemma \ref{P:ker}). We fix $k\in\N$ and consider the interval $[k\tau_a,(k+1)\tau_a]$ (whose endpoints are two consecutive zeros of $v_{a,0}$). Then, by the Sturm separation theorem (see \cite{Teschl}) we infer that the point $T_{a,k}$ is indeed the unique zero of $w_{a,0}$ lying in $[k\tau_a,(k+1)\tau_a]$. Hence $(i)$ is proved and the proof of the Lemma is complete.\qed
\end{Proof}
\begin{Remark}\label{R:Ta0}
$ \ $
\begin{itemize}
\item[$(i)$]
If $\bar{a}\in(-1,0)\cup(0,\infty)$ and $\overline{T}>0$ are such that $\overline{T}\not\in\mathcal{T}_{\bar{a}, 0}$, then the map $F(a,T)=(x_{a}(T),z_{a}(T))$ is a diffeomorphism between a neighborhood of $(\bar{a},\overline{T})$ and $F(\bar{a},\overline{T})$. This fact can be obtained with the inverse function theorem, since $F$ is smooth and its Jacobian is
\[
J_{F}(a,T)=z'_{a}(T)\frac{\partial x_{a}}{\partial a}(T)-x'_{a}(T)\frac{\partial z_{a}}{\partial a}(T)=-\frac{w_{a,0}(T)}{x_{a}(T)}\,.
\]
\item[$(ii)$]
Viceversa, if there exists a path $a\mapsto T_a$ from $(\bar{a}-\delta, \bar{a}+\delta)$ to $(0,\infty)$ such that 
\begin{equation}
\label{RL-constant}
x_a(T_a)=x_{\bar{a}}(T_{\bar{a}})\quad\text{and}\quad z_a(T_a)=z_{\bar{a}}(T_{\bar{a}})\quad\forall a\in (\bar{a}-\delta, \bar{a}+\delta)\,,
\end{equation} 
then $w_{\bar{a},0}(T_{\overline{a}})=0$. Indeed, differentiating the identities \eqref{RL-constant} with respect to $a$, we get that 
\begin{equation}
\label{RL-constant-1}
x'_a(T_a) \frac{\partial T_a}{\partial a}+\frac{\partial x_a}{\partial a}(T_a)=0\quad\text{and}\quad z'_a(T_a) \frac{\partial T_a}{\partial a}+\frac{\partial z_a}{\partial a}(T_a)=0\quad\forall a \in (\bar a-\delta, \bar a+\delta)\,. 
\end{equation}
We observe that $x'_{\bar{a}}(T_{\bar{a}})\ne 0$. Otherwise, $T_{\bar{a}}=k\tau_{\bar{a}}$ for some $k\in\mathbb{N}$. Then, by Lemma \ref{L:nod0}, $T_{\bar{a}}\not\in\mathcal{T}_{\bar{a},0}$ and we reach a contradiction thanks to part $(i)$. We can also show that $z'_{\bar{a}}(T_{\bar{a}})\ne 0$. Otherwise, by \eqref{eq:xaza}, it would be $x_{\bar{a}}^{2}(T_{\bar{a}})=\gamma_{\bar{a}}$. As a consequence, $\bar a>0$ and $T_{\bar{a}}=\frac{\tau_{\bar{a}}}{2}+k\tau_{\bar{a}}$ for some $k\in\mathbb{N}$. Again, by Lemma \ref{L:nod0}, $T_{\bar{a}}\not\in\mathcal{T}_{\bar{a},0}$ and we reach a contradiction thanks to part $(i)$. Having proved that $x'_{\bar{a}}(T_{\bar{a}})\ne 0$ and $z'_{\bar{a}}(T_{\bar{a}})\ne 0$, using \eqref{RL-constant-1}, we obtain two different expressions for $\frac{\partial T_a}{\partial a}\big|_{a=\bar{a}}$ and, equating them, we readily deduce the assertion.
\end{itemize}
\noindent
In order to explain the geometrical meaning of the above sentences, let us denote $\Sigma_{a,T}:=X_{a}([-T,T]\times \mathbb{S}^{1})$ the bounded section of the Delaunay surface parameterized by $X_{a}$, whose boundary is given by two coaxial circles of radius $x_{a}(T)$, lying in the planes $z=\pm z_{a}(T)$. 
Part $(i)$ states that if $w_{\bar{a},0}(\overline{T})\ne 0$, then there exists a neighborhood $U$ of $(\bar{a},\overline{T})$ such that in the class of surfaces $\Sigma_{{a},{T}}$ with $(a,T)\in U$, only the surface $\Sigma_{\bar{a},\overline{T}}$ has its boundary consisting of the circles defined by $x^{2}+y^{2}=\bar{R}^{2}$ and $z=\pm \bar{L}$, with $\bar{R}=x_{\bar{a}}(\overline{T})$ and $\bar{L}=z_{\bar{a}}(\overline{T})$. Viceversa, part $(ii)$ states that if there exists a path of surfaces $\Sigma_{a,T_{a}}$ with $a$ in a neighborhood of $\bar{a}$, having the same boundary, then $T_{\bar{a}}$ has to be a zero of $w_{\bar{a},0}$. Hence, roughly speaking, when $T\not\in\mathcal{T}_{a,0}$, the surface $\Sigma_{a,T}$ has a certain rigidity (i.e., nondegeneracy), up to translation in the horizontal plane. 
\end{Remark}

\noindent{\it Proof of Theorem \ref{T:Ta}.\ \ }Let $a \in (-1,0) \cup(0, \infty)$ be fixed. The inclusion 
\[
\bigcup_{j=0}^{\infty}\mathcal{T}_{a,j} \subseteq \T_a
\]
is trivial. Let us prove the opposite inclusion. Let $T \in \T_a$ and let $\f \in C^2([-T, T]\times \S^1)$ be a nontrivial even solution to $\mathfrak{L}_{a}\f = 0$ with $\f(\pm T, \cdot) =0$ on $\S^1$. Since $\f$ is periodic with respect to $\theta$, we can express it via Fourier series
\[
\f(t, \theta) =\sum_{j \in \Z} \f_j(t)\, \xi_j(\theta)\,,
\]
where
\begin{equation}\label{eq:Fourier}
\xi_j(\theta) = 
\begin{cases}
\frac{1}{{\sqrt{2\pi}}} & \text{ if }j=0\\
\frac{1}{{\sqrt{\pi}}}\cos(j\theta) & \text{ if }j>0\\
\frac{1}{{\sqrt{\pi}}}\sin(j\theta) & \text{ if }j<0
\end{cases}
\qquad \text{ and } \qquad
\f_j(t ) = \int_{-\pi}^\pi \f(t, \theta)\,\xi_j(\theta) d \theta\,.
\end{equation}
By definition, for every $j \in \Z$ it holds that $\f_j \in C^{2}([-T,T])$, $\f_j(\pm T)=0$, and $\f_j(t) = \f_j(-t)$ for any $t \in [-T, T]$. Moreover, 
\[
\mathfrak{L}_{a}\varphi = \frac{1}{2x_a^2}\sum_{j \in \Z}(\f''_j(t) - (j^2-2p_a(t))\f_j(t))\ \xi_j(\theta)\,,
\]
thus $\mathfrak{L}_a \f= 0$ implies that  for any $j \in \mathbb{Z}$, $\f_j$ is an even classical solution to
\[
\begin{cases}
\f''_j = (j^2-2p_a)\f_j & \text{ in } (-T, T)\\
\f_j(\pm T) = 0\,.
\end{cases}
\]
We observe that for any $j \in \N$, since $\f_{j}$ and $\f_{-j}$ satisfy the same equation, we get that $\f_{-j} = \eta_j\f_{j}$, for some $\eta_j \in \R$. Moreover, by linearity, for every $j \in \N$ it holds $\f_j = \lambda_jw_{a,j}$ for some $\lambda_j \in \R$, where $w_{a,j}$ are the functions defined in Lemma \ref{P:ker} and \eqref{eq:wv}.  Thus, if $T \not \in \T_{a,j}$, since $\f_j(\pm T) = 0$ the only possibility is that $\lambda_j =0$, hence we conclude that $\f_j \equiv 0$.  On the other hand, since $\f$ is nontrivial there exists $j_0 \in \N$ such that $\f_{j_0}\not\equiv0$, which implies that $T \in \T_{a,j_0}$.
As a consequence $\T_a \subseteq \bigcup_{j=0}^{\infty}\mathcal{T}_{a,j}$, and we deduce that $\T_a = \bigcup_{j=0}^{\infty}\mathcal{T}_{a,j}$. 

\smallskip

To prove $(i)-(iv)$ we begin by observing that, if $a \in (-1,0) \cup \left(0, \frac{\sqrt3-1}{2}\right]$ and $j \geq 2$ we have that
\[
j^2 - 2(1+2\gamma_a) \geq 2(1-2\gamma_a) \geq 0\,.
\]
Therefore, by Corollary \ref{L:zeridist}, $(i)$ we get that $\T_{a,j} = \emptyset$ for any $j\geq 2$. Now, $(i)-(iv)$ are immediate consequences of Lemma \ref{L:nod0}. 

Finally, $(v)$ follows from Corollary \ref{L:zeridist}, $(ii)$, and the proof is complete.  \qed

As a corollary of Theorem \ref{T:Ta} we provide the following characterization of $\mathrm{ker}(\mathfrak{L}_a)$ in ${\X}_T$. 
\begin{Corollary}\label{L:kerdesc}
Let $a\in (-1,0)\cup(0,\infty)$. Then there holds
\[
\begin{aligned}
&\dim\ker ({\mathfrak{L}_{a})\big|}_{\X_T} <\infty &  &\text{for every }T >0\\
&\ker ({\mathfrak{L}_{a})\big|}_{\X_T} = \{0\} &  &\text{for every } T>0 \text{ except for a countable set}\,. 
\end{aligned}
\]
More precisely, 
\begin{enumerate}
\item[$(i)$]
if $a \in (-1, 0)$, then $\ker ({\mathfrak{L}_{a})\big|}_{\X_T} \neq \emptyset$ if and only if $T \in \T_{a,0}$. In such case,
\[
\ker( {\mathfrak{L}_{a})\big|}_{\X_T} = \mathrm{span}\{w_{a,0}(t)\}\,.
\]
\item[$(ii)$]  if $a \in \left(0, \frac{1}{2}(\sqrt3-1)\right]$, then $\ker ({\mathfrak{L}_{a})\big|}_{\X_T} \neq \emptyset$ if and only if $T \in \T_{a,0}\cup \T_{a,1}$. In such case,
\[
\ker ({\mathfrak{L}_{a})\big|}_{\X_T} =
\begin{cases} 
\mathrm{span }\{w_{a,0}(t)\} & \quad \text{ if } T \in \T_{a,0}\\
\mathrm{span }\{w_{a,1}(t)\cos \theta, w_{a,1}(t)\sin \theta\} & \quad \text{ if } T \in \T_{a,1}
\end{cases}
\]

\end{enumerate}
\end{Corollary}
\begin{Proof} Let $a \in (-1,0)\cup(0,\infty)$ and let us define 
\[
J_a := \max\{ j \in \N \ \mid j^2 < 2(1+2\gamma_a)  \}\,.
\]
Using Corollary \eqref{L:zeridist}, $(i)$, we find that $\T_{a,j} = \emptyset$ for any $j\geq J_a$. Thus, by Theorem \ref{T:Ta} we have that $\T_a$ is a finite union of nonempty sets. Moreover $\T_{a,0}, \T_{a,1}$ are discrete sets (see Theorem \ref{T:Ta}, $(i),(iii)$) and, by Lemma \ref{L:genq}, $(ii)$ with $q = j^2-2p_a$ we infer that, if $\T_{a,j}$, $j \geq2$ is not empty, then it is a discrete set too. Therefore, $\T_a$ is a finite union of discrete sets and it is countable. Hence $\ker ({\mathfrak{L}_{a})\big|}_{\X_T} = \{0\}$ for any $T>0$ except for a countable set.

Next, let us fix $T>0$. If $T \not \in \T_a$, thanks to the previous discussion we have that $\ker ({\mathfrak{L}_{a})\big|}_{\X_T} = \{0\}$, and thus it trivially holds that $\dim\ker ({\mathfrak{L}_{a})\big|}_{\X_T} <\infty$. Hence, assume that $T \in \T_a$ and let us set
\[
\mathcal{W}^T_a := \{ w_{a,j}(t)\xi_j(\theta),  w_{a,j}(t)\xi_{-j}(\theta)  \mid j \leq J_a, T \in \T_{a,j}\},  
\]
where the functions $\xi_j$ are as in \eqref{eq:Fourier}. Notice that, as proved before, $\mathcal{W}^T_a$ is a finite set. As shown in the proof of Theorem \ref{T:Ta} , we have that $\f \in C^2([-T, T]\times \S^1)$ is a nontrivial even solution to $\mathfrak{L}_a\f=0$ with $\f(\pm T, \cdot) = 0$ on $\S^1$ if and only if $\f$ is a linear combination of the elements in $\mathcal{W}_a^T$, that is, 
\[
\ker ({\mathfrak{L}_{a})\big|}_{\X_T} = \mathrm{span }\,\mathcal{W}^T_a\,,
\]
and thus $\ker( {\mathfrak{L}_{a})\big|}_{\X_T}$ has finite dimension. Finally, $(i)$, $(ii)$ are now immediate consequences of Theorem \ref{T:Ta},$(i)-(iv)$. The proof is complete. 
\qed
\end{Proof}

\section{Proof of Theorem \ref{T:existence}}

In this Section we prove Theorem \ref{T:existence}. The proof is split in two parts: in the first one we deal with the general case (Theorem \ref{T:existence}, (i)) and, in the second one, we focus on curvatures which depend only on $z$ (Theorem  \ref{T:existence}, (ii)).

\subsection{The general case}
We fix $a\in(-1,0)\cup(0,\infty)$ and we consider surfaces of annular type which are normal graphs of a compact section of an unduloid or a nodoid with vertical axis of revolution. More precisely, fixing $T>0$ and $\alpha \in (0,1)$, we consider parametric surfaces defined by maps $X\colon[-T,T]\times\R/2\pi\Z\to\mathbb{R}^{3}$ of the form
\[
X(t,\theta)=X_{a}(t,\theta)+p\mathbf{e}_{1}+q\mathbf{e}_{2}+\varphi(t,\theta)N_{a}(t,\theta)\,,
\]
where $X_{a}$ and $N_{a}$ are defined in \eqref{Xa}--\eqref{La-def}, $(p,q)\in\R^2$ and $\varphi\in C^{2,\alpha}([-T,T]\times\R/2\pi\Z,\mathbb{R})$ with $\|\varphi\|_{C^{2,\alpha}}$ small enough, so that $X$ is a regular surface, too. Since we perturb with normal variations depending also on $\theta$, the resulting surface in not necessarily a surface of revolution. Moreover, we impose null boundary conditions at $t=\pm T$ and, since we also assume that the mean curvature function is even with respect to the vertical axis, we can impose the same condition also on $\varphi$ with respect to $t$. We aim to find solutions to problem \eqref{M=H}.
The unknown is the triple $(p,q,\varphi)\in \R\times\R\times C^{2, \alpha}([-T,T]\times\R/2\pi\Z,\mathbb{R})$.

\medskip

It is convenient to introduce the H\"older spaces
\[
\begin{aligned}
&\mathcal{X}_{T}:=\{\varphi\in C^{2,\alpha}([-T,T]\times\R/2\pi\Z)\mid\varphi(\pm T,\cdot)=0\,,~\varphi(\cdot,\theta)\text{ even }\forall\theta\}\\
&\mathcal{Y}_{T}:=\{\varphi\in C^{0,\alpha}([-T,T]\times\R/2\pi\Z)\mid\varphi(\cdot,\theta)\text{ even }\forall\theta\}\,.
\end{aligned}
\]
Clearly, $\mathcal{X}_{T}$ and $\mathcal{Y}_{T}$ are Banach spaces, equipped with their standard norms.

\begin{Lemma}\label{L:invertible}
Let $a\in (-1,0)\cup(0,\infty)$ and let $T>0$. Then the operator $\mathfrak{L}_{a}$ defined in \eqref{La-def} is a bounded linear operator from $\mathcal{X}_{T}$ into $\mathcal{Y}_{T}$. Moreover, $\mathfrak{L}_{a}$ is an isomorphism between $\mathcal{X}_{T}^{\bot}$ and $\mathcal{Y}_{T}^{\bot}$, where 
\[
\begin{aligned}
&\mathcal{X}_{T}^{\bot}:=\left\{\varphi\in\mathcal{X}_{T}\ \Big|\ \int_{\mathcal{R}_T}x_a^2\,\varphi\, v\,dt\,d\theta=0 \text{ for any }v \in \ker({\mathfrak{L}_{a})\big|}_{\mathcal{X}_T}\right\}\,,\\
&\mathcal{Y}_{T}^{\bot}:=\left\{g\in\mathcal{Y}_{T}\ \Big|\ \int_{\mathcal{R}_T}x_a^2\,g\, v\,dt\,d\theta=0 \text{ for any }v \in \ker( {\mathfrak{L}_{a})\big|}_{\mathcal{X}_T}\right\}\,,
\end{aligned}
\]
with the agreement that $\X_T^\perp = \X_T$, $\mathcal Y_T^\perp = \mathcal Y_T$ if $\ker ({\mathfrak{L}_{a})\big|}_{\mathcal{X}_T} = \{0\}$. 
\end{Lemma}

\begin{Proof}
The first statement of the Lemma is obvious. Let us prove that $\mathfrak{L}_a$ is an isomorphism. To this end, we introduce the Hilbert space $\mathrm{X}$ obtained as the completion of $\mathcal X_T^\perp$ with respect to the standard norm in $W^{1,2}(\mathcal{R}_T)$. Consider the eigenvalue problem 
\begin{equation}\label{eq:eigprob}
-\Delta u + 2q u = \lambda u\,, \quad u \in \mathrm{X}\,,
\end{equation}
where $q \geq 0$ is defined as 
\[
q := (1+2\gamma_{a}) - p_a \in C^\infty([-T, T])\,.
\]
For every $u \in \mathcal{X}_T$ it holds
\[
\int_{\mathcal{R}_T} |\nabla u|^2 dt\,d\theta\geq \int_{-\pi}^\pi \left(\ \int_{-T}^T\left|\frac{\partial u}{\partial t}\right|^2 \,dt\right) \, d \theta \geq  \frac{\pi^2}{4T^2} \int_{-\pi}^\pi \left(\ \int_{-T}^T\left|\frac{\partial u}{\partial t}\right|^2 \,dt\right) \, d \theta = \frac{\pi^2}{4T^2} \int_{\mathcal{R}_T}|u|^2 \,dt\, d\theta\,.
\]
Thus by density we get the Poincar\'e-type inequality
\[
\int_{\mathcal{R}_T} |\nabla u|^2 dt\,d\theta\geq \frac{\pi^2}{4T^2} \int_{\mathcal{R}_T}|u|^2 \,dt\, d\theta\, \quad \text{ for any }u \in \mathrm{X}\,.
\]
Thanks to previous discussion and a standard argument, we infer that there exists a strictly positive sequence of eigenvalues for problem \eqref{eq:eigprob}. 

Let now $g \in \mathcal{Y}_T$, and consider the equation
\begin{equation}\label{eq:surj}
\mathfrak{L}_{a} u = g\,,
\end{equation}
which is equivalent to 
\[
-\Delta u + q u = 2(1+2\gamma_{a}) u - x_a^2 g\,.
\]
Once again, by a standard argument we get a Fredholm alternative-type result: if $2(1+2\gamma_{a})$ is not an eigenvalue of \eqref{eq:eigprob}, then there exists a unique solution $u \in \mathrm{X}$ to \eqref{eq:surj}. Otherwise, if $2(1+2\gamma_{a})$ is an eigenvalue, denoting by $E$ the associated eigenspace, \eqref{eq:surj} admits a solution if and only if 
\begin{equation}\label{eq:perp}
\int_{\mathcal{R}_T}x_a^2\,g\,v\  dt\, d\theta= 0\, \quad \text{ for any }v \in E. 
\end{equation}
In this case, if $\tilde u \in \mathrm{X}$ is a solution, every other solution is in the form $u = \tilde u + v $ for some $v \in E$. Moreover, in both cases, by standard regularity theory, we get that $u \in \mathcal{X}_T$.

Next we notice that $2(1+2\gamma_{a})$ is an eigenvalue of \eqref{eq:eigprob} if and only if $T \in \T_a$. If $T \not \in \T_a$, then $\ker( {\mathfrak{L}_{a})\big|}_{\mathcal{X}_T} = \{0\}$, thus $\mathcal{X}_T = \mathcal{X}_T^\perp$, $\mathcal{Y}_T = \mathcal{Y}_T^\perp$, hence the proof is complete. On the other hand, if $T \in \T_a$ then 
$E = \ker ({\mathfrak{L}_{a})\big|}_{\mathcal{X}_T}$ and \eqref{eq:perp} is equivalent to $g \in \mathcal{Y}_T^\perp$. 
Let $\tilde u$ be a solution to \eqref{eq:surj}. It holds that (see the proof of Corollary \ref{L:kerdesc}),
\[
\ker ({\mathfrak{L}_{a})\big|}_{\mathcal{X}_T} = \text{span}\{ w_{a,j}(t)\xi_j(\theta),  w_{a,j}(t)\xi_{-j}(\theta)  \mid j \leq J_a, T \in \T_{a,j}\}\,,  
\]
where $\xi_j$ are as in \eqref{eq:Fourier}, and we check that 
\[
\int_{\mathcal{R}_T} x_a^2 w_{i,a}w_{j,a} \xi_{\pm i}(\theta)\xi_{\pm j}(\theta) = \left(\int_{-T}^T x_a^2 w_{i,a}w_{j,a} dt \right) \delta_{ij}\,,\ \quad \text{for any }i, j \in \N, i \leq J_a, j \leq J_a\,,
\]
where $\delta_{ij}$ is the Kronecker delta. 
Therefore
\[
u(t, \theta) := \tilde u(t, \theta) - \sum_{\underset{|\ell| \leq J_a}{\ell \in \mathbb{Z}}} \frac{\int_{\mathcal{R}_T} x_a^2\, \tilde u\, w_{a,|\ell|}\xi_\ell\, dt\, d\theta}{\int_{-T}^T x_a^2\, |w_{|\ell|,a}|^2dt} w_{a,|\ell|}(t)\xi_\ell(\theta)
\]
is also a solution to \eqref{eq:surj}, and it is the unique solution belonging to $\mathcal{X}_T^\perp$. 
Therefore $\mathfrak{L}_a$ is an invertible bounded operator between the Banach spaces $\X_T^\perp$ and $\mathcal{Y}_T^\perp$ and thus it is an isomorphism. The Lemma is proved.\qed 
\end{Proof}
We are now in position to prove Theorem \ref{T:existence}, $(i)$.
\medskip
\bigskip

\noindent{\it Proof of Theorem \ref{T:existence}, $(i)$.\ \ } 
Fix $a\in (-1,0)\cup(0,\infty)$, $T>0$ with $T \not \in \T_a$ and $(p,q) \in \R^2$. Set
\[
\F(\varepsilon,\varphi)=\mathfrak{M}(X_{a}+p\mathbf{e}_1 + q\mathbf{e}_2 + \varphi N_{a})-H_{\varepsilon}(X_{a}+p\mathbf{e}_1 + q\mathbf{e}_2+\varphi N_{a})\,,\quad(\varepsilon,\varphi)\in\R\times\mathcal{X}_{T}\,.
\]
By direct computation we see that $\f \mapsto \mathfrak{M}(X_{a}+p\mathbf{e}_1 + q\mathbf{e}_2+\varphi N_{a})$ maps $\mathcal{X_T}$ to $\mathcal{Y}_T$ and it is of class $C^1$. Moreover, using that $H_\e(x,y,z)$ is of class $C^{1}$ in $\R \times \R^{3}$ and since $(H_2)$ holds, we find that $\f \mapsto H_\e(X_{a}+p\mathbf{e}_1 + q\mathbf{e}_2+\varphi N_{a})$ maps $\mathcal{X_T}$ to $\mathcal{Y}_T$ and it is of class $C^1$. Thus we conclude that $\F \in C^1(\R\times\mathcal{X}_{T};\mathcal{Y}_{T})$. 

Since $\mathfrak{M}(X_a+p\mathbf{e}_1 + q\mathbf{e}_2) \equiv \mathfrak{M}(X_a) \equiv {H}_{0}\equiv 1$, we have
\[
\F(0,0)=0\,.
\]
In addition, since $\nabla H_0\equiv 0$, using also \eqref{La-def}, we obtain  
\[
\frac{\partial \F}{\partial\varphi}(0,0)[\xi,\psi]=\mathfrak{L}_{a}\psi\,, \quad (\xi, \psi)\in\R\times\mathcal{X}_{T}\,.
\]
Thanks to Lemma \ref{L:invertible} we can apply the Implicit function theorem and conclude that there exist $\overline\varepsilon=\overline\varepsilon(a,T,p,q)$ and a function $\varepsilon\mapsto\varphi_{\varepsilon}\in\mathcal{X}_{T}$ defined for $|\varepsilon|<\overline\varepsilon$ of class $C^{1}$, such that $\f_{0}=0$ and $\F(\varepsilon,\varphi_{\varepsilon})=0$ for every $\varepsilon\in(-\overline\varepsilon, \overline\varepsilon)$. This means that $\varphi_{\varepsilon}$ solves \eqref{M=H} and the proof is complete. 
\qed

\subsection{The case of mean curvature functions depending just on $z$}

In this Subsection we consider mappings $H_\e$ of the form 
\[
H_{\varepsilon}(x,y,z)=\widetilde{H}_\e(z) \quad \forall (x, y, z) \in \R^3\,,
\]
where $\widetilde{H}_\e\colon \R\times\R\to\R$ is a mapping in $C^{1}(\R \times \R)$ such that
\begin{equation}
\begin{gathered}\label{eq:condH}
\widetilde{H}_0(z)=1\quad\forall z\in\R\\
\widetilde{H}_\e(z)=\widetilde{H}_\e(-z)\quad\forall z\in\R\,,~\forall\varepsilon\in\R\,.
\end{gathered}
\end{equation}

We fix $a\in(-1,0)\cup(0,\infty)$ and we consider surfaces of revolution of annular type which are normal graphs of $X_a([-T, T]\times \S^1)$. More precisely, fixing $T>0$, we consider parametric surfaces defined by maps $X\colon[-T,T]\times\mathbb{S}^{1}\to\mathbb{R}^{3}$ of the form
\[
X(t,\theta)=X_{a}(t,\theta)+\varphi(t)N_{a}(t,\theta)
\]
where $\varphi\in C^{2}([-T,T],\mathbb{R})=:C^{2}[-T,T]$. Also in this case we aim to find solutions to the problem \eqref{M=H}.
Here the unknown is the triple $(p,q,\varphi)\in \R\times\R\times C^{2}[-T,T]$.

The fact that we take variations along the normal vector depending just on $t$ guarantees that $X(t,\theta)$ yields a parameterization of a surface of revolution. This choice is meaningful because we are studying a problem of prescribed mean curvature in which the mean curvature function depends just on the $z$ variable. In this case, the problem simplifies a lot because is reduced to an ODE. Moreover, we impose null boundary conditions and, since we also assume that the mean curvature function is even, we ask $\varphi$ to satisfy the same symmetry. Hence we introduce the Banach spaces
\[
\widetilde{\mathcal{X}}_{T}:=\{\varphi\in C^{2}[-T,T]\mid\varphi(\pm T)=0\,,~\varphi\text{ even}\}\quad\text{and}\quad
\widetilde{\mathcal{Y}}_{T}:=\{\varphi\in C^{0}[-T,T]\mid\varphi\text{ even}\}
\]
equipped with their standard norms. We have that:

\begin{Lemma}\label{L:opnorminvJacobi}
Let $a\in (-1,0)\cup(0,\infty)$ and let $T>0$, $T \not \in \T_{a,0}$. Then the operator $\mathfrak{L}_{a}$ defined in \eqref{La-def} is an isomorphism of Banach spaces between $\widetilde{\mathcal{X}}_{T}$ and $\widetilde{\mathcal{Y}}_{T}$. 
\end{Lemma}
\begin{Proof}
 Let us fix $a\in (-1,0)\cup(0,\infty)$ and $T>0$ with $w_{a,0}(T)\ne 0$. By Proposition \ref{L:Jacobi} we have that
\[
\mathfrak{L}_{a}\varphi=\frac{1}{2x_{a}^{2}}\left(\varphi''+2p_{a}\varphi\right)\quad\forall\varphi\in\widetilde{\mathcal{X}}_{T}\,.
\]
Then, we check that $\mathfrak{L}_{a}$ is a bounded linear operator from $\widetilde{\mathcal{X}}_{T}$ into $\widetilde{\mathcal{Y}}_{T}$. By Lemma \ref{P:ker}, if $\varphi\in\widetilde{\mathcal{X}}_{T}$ and $\mathfrak{L}_{a}\varphi=0$, then $\varphi=c_1w_{a,0} + c_2 v_{a,0}$ for some constants $c_{1}\,,~c_{2}\in\R$. Since $\varphi$ is even, then, $c_{2}=0$ and as $\varphi(\pm T)=0$, we also infer that $c_{1}=0$. Hence the operator $\mathfrak{L}_{a}$ has null kernel in $\mathcal{X}_{T}$. Moreover, for $g\in\widetilde{\mathcal{Y}}_{T}$, the function
\[
\begin{split}
\varphi(t)&=\frac{2}{1+2a}\left\{\left(\int_{0}^{t}x_a^2\,g\,w_{a,0}\,ds\right)v_{a,0}(t)-\left(\int_{0}^{t}x_a^2\, g\, v_{a,0}\,ds\right)w_{a,0}(t)\right.\\
&\left.+\frac{1}{w_{a,0}(T)} \left[-\left(\int_{0}^{T}x_{a}^{2}\,g\,w_{a,0}\,ds\right)v_{a,0}(T)+\left(\int_{0}^{T}x_a^2\, g\, v_{a,0}\,ds\right)w_{a,0}(T)\right]w_{a,0}(t)\right\}
\end{split}
\]
solves $\mathfrak{L}_{a}\varphi=g$ and belongs to $\widetilde{\mathcal{X}}_{T}$, and the conclusion plainly follows.\qed
\end{Proof}

\begin{Remark}
With a view to studying problems involving the singular limit $a\to 0$, it is important to estimate the norm of the inverse of the operator $\mathcal{L}_{a}\colon\widetilde{\mathcal{Y}}_{T}\to\widetilde{\mathcal{X}}_{T}$ with respect to small $|a|$. One can show that 
\begin{equation}
\label{norm-La-1}
\|\mathcal{L}_{a}^{-1}\|\to\infty\quad\text{as $a\to 0$\,,}
\end{equation} 
if $T$ is larger than the only positive zero $T_{0}$ of the function $w_{0}(t)=-1+t\tanh(t)$. Indeed, one has that $w_{a,0}\to w_{0}$ in $C^2_{loc}(\R)$, as $a\to 0$ (see \cite[Lemma 4.7]{CaldiroliMusso}, \cite[Lemma 2.8]{CIM}). In particular, $w_{a,0}(T_{0})\to 0$, as $a\to 0$. Then, taking $g\in\widetilde{\mathcal{Y}}_{T_{0}}$ such that 
\[
\int_{0}^{T_{0}}g(s) (1-s\tanh(s)) \sech(s)^2\,ds\neq 0\,,
\]
since $v_{a,0}(t)\to \tanh(t)$ as $a\to 0$, we have that
\[
\left|\frac{v_{a,0}(T_0)}{w_{a,0}(T_{0})} \int_{0}^{T_{0}}x_{a}^{2}(s)\,g(s)\,w_{a,0}(s)\,ds\right|\to\infty
\]
and then \eqref{norm-La-1} easily follows. 
A similar reasoning holds true also for fixed $a\in(-1,0)\cup(0,\infty)$ and $T$ varying in a neighborhood of some $T_{a,k}\in \T_{a,0}$. 
\end{Remark}

\noindent{\it Proof of Theorem \ref{T:existence}, $(ii)$.\ \ } 
It suffices to argue as in the proof of Theorem \ref{T:existence}, $(i)$, taking into account Lemma \ref{L:opnorminvJacobi}. \qed

\section{Proof of Theorem \ref{T:nonexistence}}

\noindent{\it Proof of Theorem \ref{T:nonexistence}.\ \ } 
Let $a \in (-1, 0) \cup (0, \infty)$, $(p_0, q_0) \in \R^2$ and let $H_\e$ be as in the statement. Here it is convenient to denote
\[
H_\e(x, y, z) = H(\e, x, y, z)\,.
\]

We prove $(i)$. Let $T = T_{a,0}$, where $T_{a,0}$ is as in Theorem \ref{T:Ta}, $(i)$. We point out that, by Lemma \ref{L:nod0}, $(i)$,  $T_{a,0}$ is the first positive zero of $w_{a,0}$ where $w_{a,0}$ is defined in Lemma \ref{P:ker}.
  
 Assume by contradiction that there exist $\bar\varepsilon>0$, a $C^0$-function $\e \mapsto (p_\e, q_\e)$ and a $C^1$-function $\varepsilon\mapsto\varphi_{\varepsilon}\in C^2([-T_{a,0},T_{a,0}]\times\S^1)$, defined for $|\varepsilon|<\bar\varepsilon$, such that $\left.(p_\e, q_\e)\right|_{\e = 0} = (p_0, q_0)$, $\varphi_{0}=0$ and for all $\varepsilon\in(-\bar\varepsilon, \bar\varepsilon)$ the triple $(p_\e, q_\e,\varphi_\e)$  satisfies 
\begin{equation}\label{eq1:teononexist}
\begin{cases}
\mathfrak{M}(X_{a}+p_\e\mathbf{e}_{1}+q_\e\mathbf{e}_{2}+\varphi_{\varepsilon}N_{a})=H(\e, X_{a}+p_\e\mathbf{e}_{1}+q_\e\mathbf{e}_{2}+\varphi_{\varepsilon}N_{a})\text{ in }[-T_{a,0},T_{a,0}]\times\mathbb{S}^{1}\\
\varphi_{\varepsilon}(\pm T_{a,0},\cdot)=0\,.
\end{cases}
\end{equation}
Since $\mathfrak{M}(X_a +  p_\e \mathbf{e}_1 + q_\e \mathbf{e}_2 + \f_\e N_a) = \mathfrak{M}(X_{a}+\veps N_{a})$, 
$\mathfrak{M}(X_{a})\equiv 1$ and $H(0,x,y,z)\equiv 1$, we can write the equation in \eqref{eq1:teononexist} as 
\begin{equation}\label{eq2:teononexist}
\mathfrak{M}(X_{a}+\veps N_{a})-\mathfrak{M}(X_{a})-\mathfrak{L}_{a}\veps + \mathfrak{L}_{a}\veps=H(\e, X_{a}+ p_\e \mathbf{e}_1 + q_\e \mathbf{e}_2+\veps N_{a})-H(0, X_{a}+\veps N_{a}).
\end{equation}
Clearly, since $\mathfrak{L}_{a} w_{a,0}=0$, $w_{a,0}(\pm T_{a,0})=0$, $\veps(\pm T_{a,0},\cdot)=0$ and since $\varphi_\e$ is $2\pi$-periodic with respect to the variable $\theta$, integrating by parts it follows that
\[
\int_{\mathcal{R}_{T_{a,0}}} x_a^2\,w_{a,0}\,{\mathfrak{L}_{a}}\veps\,dt\,d\theta=\int_{\mathcal{R}_{T_{a,0}}}x_a^2\,\f_\e\,{\mathfrak{L}_{a}}w_{a,0}\,dt\,d\theta=0\,.
\]
Hence, multiplying each side of \eqref{eq2:teononexist} by $x_a^2w_{a,0}$ and integrating on $[-T_{a,0},T_{a,0}]\times\S^1$ we get that for all $\varepsilon\in(-\bar\varepsilon, \bar\varepsilon)$ it holds
\begin{equation}\label{eq3:teononexist}
\begin{aligned}
\int_{\mathcal{R}_{T_{a,0}}} &x_a^2\, w_{a,0}\,\left(\mathfrak{M}(X_{a}+\veps N_{a})-\mathfrak{M}(X_{a})-\mathfrak{L}_{a}\veps\right)\,dt\,d\theta\\
=&\int_{\mathcal{R}_{T_{a,0}}} x_a^2\,w_{a,0}\,\left(H(\e, X_{a}+ p_\e \mathbf{e}_1 + q_\e \mathbf{e}_2+\veps N_{a})-H(0, X_{a}+\veps N_{a})\right)\,dt\,d\theta\,.
\end{aligned}
\end{equation}

For the left-hand side of \eqref{eq3:teononexist} we first observe that since $\e\mapsto\veps$ is a $C^1$-mapping and as $\varphi_0=0$, then it follows that
\begin{equation}\label{eq5:teonoexist}
\|\varphi_\e\|_{C^2([-T_{a,0},T_{a,0}]\times\S^1)}\leq C_1 |\e|,
\end{equation}
for all sufficiently small $\e$, where $C_1$ is a positive constant independent of $\e$. Now, since the map $\varphi\mapsto \mathfrak{M}(X_{a}+\varphi N_{a})$ is well defined and differentiable, as a function from $\{\varphi \in C^2([-T_{a,0},T_{a,0}]\times\S^1)\mid \ \|\varphi\|_{C^2}<\delta\}$ to $C^0([-T_{a,0},T_{a,0}]\times\S^1)$, where $\delta>0$ is small number depending only on $a$, then, by definition of $\mathfrak{L}_{a}$ and thanks to \eqref{eq5:teonoexist} we infer that for all sufficiently small $\e$
\begin{equation}\label{eq6add:teonoexist}
\|\mathfrak{M}(X_{a}+\veps N_{a})-\mathfrak{M}(X_{a})-\mathfrak{L}_{a}\veps\|_{C^0([-T_{a,0}, T_{a,0}]\times\S^1)}=o(\|\varphi_\e\|_{C^2([-T_{a,0},T_{a,0}]\times\S^1)}).
\end{equation}
Hence, combining \eqref{eq6add:teonoexist} and \eqref{eq5:teonoexist} we conclude that, as $\e\to 0$, 
\begin{equation}\label{eq6:teonoexist}
\int_{\mathcal{R}_{T_{a,0}}} x_a^2\, w_{a,0}\,\left(\mathfrak{M}(X_{a}+\veps N_{a})-\mathfrak{M}(X_{a})-\mathfrak{L}_{a}\veps\right)\,dt\,d\theta=o(\e)\,.
\end{equation}

We now analyze the right-hand side of \eqref{eq3:teononexist}. We begin observing that since the map $\e\mapsto H(\e,\cdot)$ is differentiable, then, for any $(x,y,z)\in\R^3$, $\e\in \R$ we have
$$ H(\e,x,y,z)- H(0,x,y,z)=\frac{\partial  H}{\partial\e}(\xi,x,y,z)\e,$$
for some $\xi=\xi(\e,x,y,z)\in\R$ such that $|\xi|\leq |\e|$. From this, and using again that $H(0, \cdot) \equiv 1$, we readily infer that
\[
\begin{aligned}
\int_{\mathcal{R}_{T_{a,0}}} x_a^2\,w_{a,0}\,\left(H(\e, X_{a}+ p_\e \mathbf{e}_1 + q_\e \mathbf{e}_2+\veps N_{a})-H(0, X_{a}+\veps N_{a})\right)\,dt\,d\theta&\\
 = \e \int_{\mathcal{R}_{T_{a,0}}} x_a^2\,w_{a,0}\, \frac{\partial H}{\partial\e}\left(\xi,X_{a}+ p_\e \mathbf{e}_1 + q_\e \mathbf{e}_2+\veps N_{a}\right)\,dt\,d\theta.&
\end{aligned}
\]
Now, since $\e \mapsto (p_\e, q_\e)$ is continuous and using \eqref{eq5:teonoexist}, we deduce that
\[
X_{a}+ p_\e \mathbf{e}_1 + q_\e \mathbf{e}_2+\veps N_{a} \to X_{a}+ p_0 \mathbf{e}_1 + q_0 \mathbf{e}_2\,,\quad \text{ as }\e\to 0\,,\quad \text{ uniformly on }[-T_{a,0},T_{a,0}]\times\S^1\,.
\]
Moreover, since $\xi=\xi(\e,X_{a}+p_\e \mathbf{e}_1+q_\e \mathbf{e}_2 + \veps N_{a})$ is such that $|\xi|\leq|\e|$, for all $\e\in(-\bar{\e},\bar{\e})$ and  $(t,\theta)\in [-T_{a,0},T_{a,0}]\times\S^1$, then, exploiting the continuity of the map $(\e,x,y,z)\mapsto \frac{\partial H}{\partial\e}(\e,x,y,z)$ we infer that
\[
\frac{\partial  H}{\partial\e}\left(\xi,X_{a}+ p_\e \mathbf{e}_1 + q_\e \mathbf{e}_2+\veps N_{a}\right)\to  \frac{\partial  H}{\partial\e}\left(0,X_a+ p_0 \mathbf{e}_1 + q_0 \mathbf{e}_2\right), \ \hbox{as $\e\to 0$}
\]
uniformly on on $[-T_{a,0},T_{a,0}]\times\S^1$. As a consequence we obtain that, as $\e \to 0$, 
\begin{equation}\label{eq7:teonoexist}
\begin{aligned}
\int_{\mathcal{R}_{T_{a,0}}} x_a^2\,w_{a,0}\,\left(H(\e, X_{a}+ p_\e \mathbf{e}_1 + q_\e \mathbf{e}_2+\veps N_{a})-H(0,X_{a}+\veps N_{a})\right)\,dt\,d\theta&\\
=\e\int_{\mathcal{R}_{T_{a,0}}} x_a^2\,w_{a,0}\, \frac{\partial H}{\partial\e}\left(0,X_{a}+ p_0 \mathbf{e}_1 + q_0 \mathbf{e}_2\right)\,dt\,d\theta + o(\e)&\,.
\end{aligned}
\end{equation}
Plugging \eqref{eq6:teonoexist}, \eqref{eq7:teonoexist} in \eqref{eq3:teononexist}, dividing by $\e$ and taking the limit as $\e \to 0$ we get
\[
\int_{\mathcal{R}_{T_{a,0}}} x_a^2\,w_{a,0} \,\left.\frac{\partial H}{\partial\e}\left(X_{a}+ p_0 \mathbf{e}_1 + q_0 \mathbf{e}_2\right)\right|_{\e = 0} dt\,d\theta= 0\,. 
\]
Now, thanks to the assumption \eqref{eq:deHzero}, since by construction the function $w_{a,0}$ has constant sign on $(-T_{a,0},T_{a,0})$, and since $x_a>0$ on $\R$, the integral on the left hand side is non null. Thus we readily get a contradiction and the proof of $(i)$ is complete.

The proof of $(ii)$ is the same, with slight adjustments. \qed

\begin{Remark}
A sufficient condition for
\[
\int_{\mathcal{R}_{T_{a,k}}}x_{a}^{2}w_{a,0}\,\left.\frac{\partial H_\varepsilon}{\partial \varepsilon}(X_a + p_0 \mathbf{e}_1 +q_0 \mathbf{e}_2 )\right|_{\varepsilon=0}\,dt\,d\theta\ne 0
\]
to hold is that $\frac{ \widetilde H}{\partial\e}\left(0,X_{a}+ p_0 \mathbf{e}_1 + q_0 \mathbf{e}_2\right)$ has the same sign of $w_{a,0}$ in its first $k$ nodal regions, i.e. if for all $(t,\theta)\in(-T_{a,k},T_{a,k})\setminus\{T_{a,j}\}_{j=1}^{k-1}\times\S^1$ it holds
$$w_{a,0} \frac{\partial H}{\partial\e}\left(0,X_{a}+ p_0 \mathbf{e}_1 + q_0 \mathbf{e}_2\right)>0.$$
 Specularly, the same conclusion holds whenever $\frac{\partial H}{\partial\e}\left(0,X_{a}+ p_0 \mathbf{e}_1 + q_0 \mathbf{e}_2\right)$ has the opposite sign of $w_{a,0}$ in its first $k$ nodal regions.
\end{Remark}

\section{Proof of Theorem \ref{T:existence2}}

In this section we prove Theorem \ref{T:existence2}. We begin with a preliminary result in which we compute the partial derivatives of the Melnikov-type function defined in \eqref{M}. 

\begin{Lemma}\label{L:Melnikov}
For any $a\in(-1,0)\cup(0,\infty)$ and $T>0$ it holds $M \in  C^1(\R^2)$. In particular
\[
\begin{gathered}
\frac{\partial M }{\partial p}(p,q)=-\int_{\mathcal{R}_T}x_a^2\,{\widetilde H}(X_{a}+p\mathbf{e}_{1}+q\mathbf{e}_{2})\, w_{a,1}\cos\theta\,dt\,d\theta\,,\\
\frac{\partial M}{\partial q}(p,q)=-\int_{\mathcal{R}_T}x_a^2\,{\widetilde H}(X_{a}+p\mathbf{e}_{1}+q\mathbf{e}_{2})\, w_{a,1}\sin\theta\,dt\,d\theta\,.
\end{gathered}
\]
If in addition ${\widetilde H} \in C^1(\R^3)$ then  $M \in  C^2(\R^2)$.
\end{Lemma}
\begin{Proof}
First we point out that the vector field $Q$, defined in \eqref{eq:Q}, is such that $\mathrm{div\, }Q = {\widetilde H}$ and $Q\cdot \mathbf{e}_3 \equiv 0$.
Moreover, the Melnikov-type function $M$ satisfies 
\begin{equation}\label{eq:melvol}
M(p,q) = \mathcal{V}_{T, Q}(X_{a}+p\mathbf{e}_{1}+q\mathbf{e}_{2})\,,
\end{equation}
where $\V_{T,Q}$ is the functional defined in \eqref{H-vol-def} (see Appendix \ref{A:areavol} for further details).  

Let us show the formula for $\frac{\partial M}{\partial p}$. Fixing $a\in(-1,0)\cup(0,\infty)$ and taking \eqref{eq:melvol} into account, we apply Lemma \ref{L:first-variation-volume} with the curve $s\mapsto X(s):=X_{a}+(s+p)\mathbf{e}_{1}+q\mathbf{e}_{2}$. Hence, since  $\dot{X}=\mathbf{e}_{1}$, $\dot{X}\cdot X_{t}\wedge X_{\theta}=x_{a}z'_{a}\cos\theta=x_{a}^{2}\,w_{a,1}\cos\theta$ and $X_{\theta}\wedge\dot{X}=-x_{a}\cos(\theta)\mathbf{e}_{3}$, it follows that
\[
\begin{aligned}
\frac{\partial M}{\partial p}(p,q) = &-\int_{\mathcal{R}_{T}}x_a^2\,{\widetilde H}(X_{a}+p\mathbf{e}_{1}+q\mathbf{e}_{2})\,w_{a,1}\cos\theta\,dt\,d\theta\\
&+\int_{-\pi}^{\pi}(Q(X_{a}+p\mathbf{e}_{1}+q\mathbf{e}_{2})\cdot \mathbf{e}_3\ x_a \cos \theta)(T,\theta)\,d\theta\\
&-\int_{-\pi}^{\pi}(Q(X_{a}+p\mathbf{e}_{1}+q\mathbf{e}_{2})\cdot \mathbf{e}_3\ x_a \cos \theta)(-T,\theta)\,d\theta\,,
\end{aligned}
\]
and, as $Q\cdot \mathbf{e}_3 \equiv 0$, we readily get the desired result. Arguing in a similar way, we obtain the expression for $\frac{\partial M}{\partial q}$. Finally, exploiting the Lebesgue dominated convergence Theorem, we obtain that both partial derivatives are continuous, and thus $M \in C^1(\R^2)$.

If we further assume that ${\widetilde H} \in C^1(\R^3)$, then by a standard argument involving the derivatives of parameter-depending integrals we obtain that $M \in C^2(\R^2)$, and the proof is complete.
\qed
\end{Proof}

\noindent{\it Proof of Theorem \ref{T:existence2}, $(i)$.\ \ } 
Let $\varepsilon\mapsto(p_{\varepsilon},q_{\varepsilon})$, $\varepsilon\mapsto\varphi_{\varepsilon}$ be as in the statement. Arguing exactly as in the proof of Theorem \ref{T:nonexistence}, $(i)$, where now we test 
\[
\begin{cases}
\mathfrak{M}(X_{a}+p_\e\mathbf{e}_{1}+q_\e\mathbf{e}_{2}+\varphi_{\varepsilon}N_{a})=H(\e, X_{a}+p_\e\mathbf{e}_{1}+q_\e\mathbf{e}_{2}+\varphi_{\varepsilon}N_{a})\text{ in }[-T,T]\times\mathbb{S}^{1}\\
\varphi_{\varepsilon}(\pm T,\cdot)=0\,
\end{cases}
\]
with, respectively, $x_a^2\, w_{a,1}\cos\theta$  and $x_a^2\, w_{a,1}\sin\theta$, we obtain
\[
\begin{aligned}
&\int_{\mathcal{R}_T}x_a^2\,{\widetilde H}(X_{a}+p_0\mathbf{e}_{1}+q_0\mathbf{e}_{2})\, w_{a,1}\cos\theta\,dt\,d\theta = 0\,,\\
&\int_{\mathcal{R}_T}x_a^2\,{\widetilde H}(X_{a}+p_0\mathbf{e}_{1}+q_0\mathbf{e}_{2})\, w_{a,1}\sin\theta\,dt\,d\theta = 0\,.
\end{aligned}
\]
Then, applying Lemma \ref{L:Melnikov}, we get
\[
\frac{\partial M}{\partial p}(p_0, q_0) = \frac{\partial M}{\partial q}(p_0, q_0) = 0\,.
\]
The proof is complete.  \qed

We now prove Theorem \ref{T:existence2}, $(ii)$. We begin with a preliminary Lemma. 

\begin{Lemma}\label{L:fdr}
Let $a \in \left(0, \frac{1}{2}(\sqrt3-1)\right]$, $T \in \T_{a,1}$, let $H_\e \in C^1(\R \times \R^3)$ satisfy $(H_1)$, $(H_2)$, $(H_3)$, and fix $R>0$.  Then there exist $\bar\varepsilon=\bar\varepsilon(a,T,R)$ and mappings $(\varepsilon,p,q)\mapsto\varphi_{\varepsilon,p,q}\in\mathcal{X}_{T}$, $(\varepsilon,p,q)\mapsto\lambda_{i,\varepsilon}(p,q)\in\R$ ($i=1,2$), defined for $|\varepsilon|<\bar\varepsilon$ and $(p,q)\in B_{R}$, of class $C^{1}$, such that $\varphi_{0,p,q}=0$, $\lambda_{i,0}(p,q)=0$ ($i=1,2$), and for every $\varepsilon\in(-\bar\varepsilon, \bar\varepsilon)$ it holds
\begin{eqnarray*}
&&\mathfrak{M}(X_{a}+p\mathbf{e}_{1}+q\mathbf{e}_{2}+\varphi_{\varepsilon,p,q} N_{a})-H_{\varepsilon}(X_{a}+p\mathbf{e}_{1}+q\mathbf{e}_{2}+\varphi_{\varepsilon,p,q} N_{a})\\
&=&\lambda_{1,\varepsilon}(p,q)w_{a,1}\cos\theta+\lambda_{2,\varepsilon}(p,q)w_{a,1}\sin\theta\,.
\end{eqnarray*}
\end{Lemma}

\begin{Proof}
Let $a$, $T$, $H_\e$ and $R$ be as in the statement. Let $(p, q) \in B_R$ be fixed. For simplicity, we set 
\[
U_{p,q} := X_{a}+p\mathbf{e}_{1}+q\mathbf{e}_{2}\,.
\]
Let us consider the map $\mathcal F=(\F_1, \F_2): \R \times \R^2 \times  \mathcal{X}_T \times \R\times \R \to C^{0,\alpha}([-T, T]\times \S^1) \times \R\times \R\,$, defined by 
\[
\begin{gathered}
\mathcal{F}_{1}(\varepsilon,p,q; \f, \lambda_1, \lambda_2) := \mathfrak{M}(U_{p,q}+\varphi N_{a})-H_{\varepsilon}(U_{p,q}+\varphi N_{a})-\lambda_{1}w_{a,1}\cos\theta-\lambda_2w_{a,1}\sin\theta\,,\\
\mathcal{F}_{2}(\f) := \left(\ \int_{\mathcal{R}_T}x_a^2 \, \f\, w_{a,1}\cos \theta\, dt\, d\theta\,,\ \int_{\mathcal{R}_T}x_a^2 \, \f\, w_{a,1}\sin \theta\, dt\, d\theta \right)\,.
\end{gathered}
\]
We aim to apply the implicit function theorem to $\mathcal F$. To this end, we first we point out that, as in the proof of Theorem \ref{T:existence}, $(i)$, we see that
\[
\quad \mathcal F \in C^1(\R \times \R^2 \times \mathcal{X}_T \times \R \times \R;\mathcal{Y}_T \times \R\times \R)\,.
\]
Moreover, since $\mathfrak{M}(U_{p,q}) = \mathfrak{M}(X_{a}) = 1$, we have
\[
\mathcal{F}(0, p,q; 0, 0, 0) = (\mathfrak{M}(X_{a}) - H_0, 0, 0 ) = (0, 0, 0)\,. 
\]

Next, by direct computation, arguing as in the proof of Theorem \ref{T:existence}, $(i)$, and taking into account that $\nabla H_0 = 0$, we check that
\[
\frac{\partial \mathcal F}{\partial (\f, \lambda_1, \lambda_2)}(0, p,q; 0,0,0)[\psi, \mu_1, \mu_2] = \mathfrak{G}(\psi, \mu_1 , \mu_2)\,, \quad \forall\, \psi \in \X_T, \mu_1, \mu_2 \in \R\,,
\]
where $\mathfrak G : \mathcal{X}_{T}\times \R\times \R \to \mathcal{Y}_T \times \R\times \R$ is the map given by 
\[
\begin{gathered}
\mathfrak{G}(\psi, \mu_1 , \mu_2) =  \left(\!\mathfrak{L}_a \psi - \mu_1 w_{a,1}\cos\theta-\mu_2w_{a,1}\sin\theta, \int_{\mathcal{R}_T}\!\!x_a^2 \, \psi\, w_{a,1}\,\cos \theta\, dt\, d\theta,\ \int_{\mathcal{R}_T}\!\!x_a^2 \, \psi\, w_{a,1}\,\sin \theta\, dt\, d\theta\! \right).
\end{gathered}
\]
We claim that $\mathfrak G$ is invertible. 

To prove the claim, we first recall that, since $a \in \left(0, \frac{1}{2}(\sqrt3-1)\right]$ and $T \in {\T}_{a,1}$, by Corollary \ref{L:kerdesc} and Lemma \ref{L:invertible} we have that 
\[
\ker({\mathfrak{L}_{a})\big|}_{\X_T} = \text{span}\{w_{a,1} \cos\theta\,,\ w_{a,1} \sin\theta\}\,,
\]
and that $\mathfrak{L}_a \colon \X_T^\perp \to \mathcal{Y}_T^\perp$ is an isomorphism, where  
\[
\begin{aligned}
&\mathcal{X}_{T}^{\bot}:=\left\{\varphi\in\mathcal{X}_{T}\ \Big|\ \int_{\mathcal{R}_T}x_a^2\,\f\, w_{a,1} \cos\theta\,dt\,d\theta= \int_{\mathcal{R}_T}x_a^2\, \f\, w_{a,1} \sin\theta\,dt\,d\theta=0 \right\}\,,\\
&\mathcal{Y}_{T}^{\bot}:=\left\{g\in\mathcal{Y}_{T}\ \Big| \ \int_{\mathcal{R}_T}x_a^2\,g\, w_{a,1}\cos\theta\,dt\,d\theta= \int_{\mathcal{R}_T}x_a^2\,g\,w_{a,1}\sin\theta\,dt\,d\theta=0 \right\}\,.
\end{aligned}
\]

Let us begin by showing that $\mathfrak G$ is injective. Assume that there exists $(\psi, \mu_1 , \mu_2) \in  \mathcal{X}_{T}\times \R\times \R$ such that  $\mathfrak{G}(\psi, \mu_1 , \mu_2) = 0$, that is
\begin{align}
&\mathfrak{L}_a \psi =\mu_1 w_{a,1}\cos\theta+\mu_2w_{a,1}\sin\theta\,,\label{eq:inj1}\\
& \int_{\mathcal{R}_T}x_a^2 \, \psi\, w_{a,1}\cos \theta\, dt\, d\theta =0\,, \quad   \int_{\mathcal{R}_T}x_a^2 \, \psi\, w_{a,1}\sin \theta\, dt\, d\theta = 0\,.\label{eq:inj2}
\end{align}
In particular, \eqref{eq:inj2} implies that $\psi \in \mathcal{X}_T^\perp$, and thus we have that $\mathfrak{L}_a \psi  \in \mathcal{Y_T^\perp}$. Then, in view of \eqref{eq:inj1}, we deduce that $\mu_1= \mu_2 = 0$. As a consequence, $\mathfrak{L}_a \psi = 0$ and $\psi \in \ker({\mathfrak{L}_{a})\big|}_{\X_T}$. Hence, since 
$\psi \in \mathcal{X}_T^\perp$, we conclude that $\psi \equiv 0$ and $\mathfrak{G}$ is injective. 

Let us prove the surjectivity of $\mathfrak{G}$. Let  us fix $(g, \nu_1, \nu_2) \in  \mathcal{Y}_T \times \R\times \R$. We want to solve 
\[
\mathfrak{G}(\psi, \mu_1 , \mu_2)  = (g, \nu_1, \nu_2 )\,.
\]
First of all, we take $\mu_1, \mu_2$ as 
\[
\mu_1 = - \frac{\int_{\mathcal{R}_T}x_a^2 \, g \, w_{a,1} \cos \theta\, dt\, d\theta }{\int_{\mathcal{R}_T}x_a^2\, |w_{a,1}|^2 |\cos \theta|^2\, dt\, d\theta} \,, \quad \mu_2 = - \frac{\int_{\mathcal{R}_T}x_a^2 \,g \, w_{a,1} \sin \theta\, dt\, d\theta }{\int_{\mathcal{R}_T}x_a^2\, |w_{a,1}|^2 |\sin \theta|^2\, dt\, d\theta}\,.
\]
Thanks to this choice, and using that $\int_{\mathcal{R}_T} x_a^2\, |w_{a,1}|^2 \cos \theta\, \sin \theta\, dt\, d\theta =0$,  we have 
\[
\tilde g := g + \mu_1 w_{a,1}\cos\theta+\mu_2w_{a,1}\sin\theta  \in \mathcal{Y}_T^\perp\,.
\]
Thus, by Lemma \ref{L:invertible} there exists a unique $\tilde  \psi  \in \mathcal{X}_T^\perp$ such that $\mathfrak{L}_a \tilde \psi  = \tilde g$. Then, defining
\[
\psi := \nu_1 w_{a,1}\cos\theta+\nu_2w_{a,1}\sin\theta + \tilde \psi\,,
\]
we see that $(\psi, \mu_1, \mu_2)$ is the solution to $\mathfrak{G}(\psi, \mu_1 , \mu_2)  = (g, \nu_1, \nu_2 )$. Hence $\mathfrak G $ is also surjective, hence a bijection. 

Applying the implicit function theorem to $\mathcal{F}$, for any fixed $(p, q) \in B_R$, and thanks to a standard compactness argument, we infer that there exist $ \bar \e = \bar\e(a, T, R) >0$
and uniquely determined ${ C}^1$ functions 
\[
\f: \left[
\begin{aligned}
&(-\bar\e,\bar\e) \times B_R \to \mathcal{X}_T \\
&(\e, p, q) \mapsto \f_{\e,p,q} 
\end{aligned}
\right.
\qquad 
\lambda_1: \left[
\begin{aligned}
& (-\bar\e,\bar\e) \times B_R \to \R \\
&(\e,p, q) \mapsto \lambda_{1,\e}(p,q)
\end{aligned}\right.
\qquad 
\lambda_2: \left[
\begin{aligned}
&(-\bar\e,\bar\e) \times B_R \to \R\\
&(\e, p, q) \mapsto \lambda_{2,\e}(p,q)
\end{aligned}\right.
\]
such that for all $(p,q)\in B_R$, $\e\in (-\bar\e, \bar\e)$
\begin{equation}\label{eq:implsol}
\f_{0, p, q} \equiv 0, \qquad \lambda_{1,0}(p, q) = 0\, , \qquad \lambda_{2,0}(p, q) = 0\, ,  \qquad\mathcal{F}(\e, p, q; \f_{\e,p,q}, \lambda_{1,\e}(p,q), \lambda_{2,\e}(p,q)) = 0\, .
\end{equation}
In particular, \eqref{eq:implsol} implies that 
\[
\begin{gathered}
\mathfrak{M}(U_{p,q}+\varphi_{\e, p, q} N_{a})-H_{\varepsilon}(U_{p,q}+\varphi_{\e, p, q} N_{a})= \lambda_{1,\e}(p,q)w_{a,1}\cos\theta+\lambda_{2, \e}(p,1)w_{a,1}\sin\theta\,,
\end{gathered}
\]
and the Lemma is proved.  \qed
\end{Proof}

\noindent{\it Proof of Theorem {\ref{T:existence2}}, $(ii)$.\ \ } Let $a$, $(p_0, q_0)$, $T$ and $H_\e$ be as in the statement. Let us fix $R >0$ such that $(p_0, q_0) \in B_R$. Let $\bar \e >0$ and $\lambda_{1,\e}$, $\lambda_{2,\e}$, $\f_{\e,p,q}$ be the mappings given by Lemma \ref{L:fdr}. Once again, for simplicity we denote $
U_{p,q} = X_a + p \mathbf{e}_1 + q \mathbf{e}_2$.

By Lemma \ref{L:fdr}, it holds
\begin{equation}\label{eq:c2ex1}
\mathfrak{M}(U_{p,q} + \f_{\e,p,q}N_a)-H_{\varepsilon}(U_{p,q}+ \f_{\e,p,q}N_a) =\lambda_{1,\varepsilon}(p,q)w_{a,1}\cos\theta+\lambda_{2,\varepsilon}(p,q)w_{a,1}\sin\theta\,.
\end{equation}
Testing \eqref{eq:c2ex1}, respectively,  with $x_a^2 w_{a,1} \cos \theta$ and $x_a^2 w_{a,1} \sin \theta$ we obtain
\[
\begin{aligned}
C_0\lambda_{1,\e}(p,q) = \int_{\mathcal{R}_T}x_a^2 (\mathfrak{M}(U_{p,q} + \f_{\e,p,q}N_a)-H_{\varepsilon}(U_{p,q} + \f_{\e,p,q}N_a) ) w_{a,1} \cos \theta \, dt\, d \theta\,,\\
C_0\lambda_{2,\e}(p,q) = \int_{\mathcal{R}_T}x_a^2 (\mathfrak{M}(U_{p,q} + \f_{\e,p,q}N_a)-H_{\varepsilon}(U_{p,q} + \f_{\e,p,q}N_a) ) w_{a,1}\sin \theta \, dt\, d \theta\,,
\end{aligned}
\]
where $C_0 = \pi\int_{-T}^T|w_{a,1}|^2\,dt$. 

Let us now define $\F : [-\overline \e, \overline \e ]\times B_R \to \R^2$ as
\[
\F(\e, p, q ) = (\e^{-1}C_0\lambda_{1,\e}(p,q), \e^{-1}C_0\lambda_{2,\e}(p,q))\,.
\]
Our aim is to apply the implicit function theorem to $\F$ in $(0, p_0, q_0)$.

First of all, we notice that, by construction, it holds that $\F \in C^1([-\overline \e, \overline \e] \setminus\{0\} \times B_R)$. 
Moreover, arguing as in the proof of Theorem \ref{T:nonexistence} and Theorem \ref{T:existence2}, $(i)$, we find that, for $\e$ sufficiently small,
\[
\F(\e, p, q) = \left( \frac{\partial M }{\partial p}(p, q),  \frac{\partial M}{\partial q}(p, q) \right) + \mathfrak{R_1}(\e, p, q)\,,
\]
where $\mathfrak{R_1}:  [-\overline \e, \overline \e ]\times B_R \to \R^2$ is of class $C^0$ and such that
\[
\sup_{(p,q) \in B_R}\mathfrak{R}_1(\e, p, q) \to 0 \quad \text{ as }\e \to 0\,.
\] 
Therefore $\F \in C^0([-\overline \e, \overline \e] \times B_R)$ and, since $(p_0, q_0)$ is a critical point of $M$, it holds
\[
\F(0, p_0, q_0 ) = 0\,.
\]

Let us compute $\frac{\partial \F_1}{\partial p }(\e, p, q )$, where $\F_1$ is the first component of $\F$. 
By definition we have 
\begin{equation}\label{eq:ex11}
\frac{\partial \F_1}{\partial p }(\e, p, q ) = \e^{-1}C_0\frac{\partial}{\partial p }\lambda_{1,\e}(p,q) \,.
\end{equation}
Next, arguing again as in the proof of Theorem \ref{T:nonexistence} and Theorem \ref{T:existence2}, $(i)$, we find
\[
\begin{aligned}
C_0\lambda_{\e,1}(p,q) &= \int_{\mathcal{R}_T}x_a^2 (\mathfrak{M}(X_a + \f_{\e,p,q}N_a)-1 ) w_{a,1} \cos \theta \, dt\, d \theta\\
&- \e\int_{\mathcal{R}_T}x_a^2 (\tilde{H}(U_{p,q} + \f_{\e,p,q}N_a)- \tilde{H}(U_{p,q}) ) w_{a,1} \cos \theta \, dt\, d \theta + \e \frac{\partial M}{\partial p }(p,q)\,.
\end{aligned}
\]
Hence
\begin{equation}\label{eq:ex12}
\begin{aligned}
C_0\frac{\partial}{\partial p }\lambda_{1,\e}(p,q) &= \int_{\mathcal{R}_T}x_a^2 \left( \frac{\partial}{\partial p }\mathfrak{M}(X_a + \f_{\e,p,q}N_a)\right) w_{a,1} \cos \theta \, dt\, d \theta\\
&- \e\int_{\mathcal{R}_T}x_a^2 \frac{\partial}{\partial p }(\tilde{H}_{\varepsilon}(U_{p,q} + \f_{\e,p,q}N_a)- \tilde{H}_{\varepsilon}(U_{p,q}) ) w_{a,1} \cos \theta \, dt\, d \theta + \e \frac{\partial^2 M}{\partial p^2 }(p,q)\,.
\end{aligned}
\end{equation}

Now, since $\tilde H \in C^2(\R^3)$, by the mean value theorem we get
\[
\begin{aligned}
\frac{\partial}{\partial p }(\tilde{H}_{\varepsilon}(U_{p,q} &+ \f_{\e,p,q}N_a)- \tilde{H}_{\varepsilon}(U_{p,q}) ) \\
& = \nabla \tilde H(U_{p,q} + \f_{\e,p,q}N_a) \cdot N_a \frac{\partial \f_{\e,p,q}}{\partial p} + \nabla\frac{\partial \tilde H}{\partial x}(X_a + p \mathbf{e}_1 + q \mathbf{e}_2 + \xi \f_{\e, p, q}N_a) \cdot N_a \f_{\e,p,q}\,,
\end{aligned}
\]
for some $\xi  \in [0,1]$.  Thus
\[
\left\|\frac{\partial}{\partial p }(\tilde{H}_{\varepsilon}(U_{p,q} + \f_{\e,p,q}N_a)- \tilde{H}_{\varepsilon}(U_{p,q}) ) \right\|_{C^0(\mathcal{R}_T)} \leq C\left( \left\| \frac{\partial \f_{\e,p,q}}{\partial p}\right \|_{C^2(\mathcal{R}_T)} + \| \f_{\e,p,q}\|_{C^2(\mathcal{R}_T)}\right)
\]
where the constant $C>0$ is independent of $(p,q)$ and $\e$. Moreover, since the map
\[
\f : (\e, p, q ) \in [-\overline \e, \overline \e] \times B_R \mapsto \f_{\e, p, q } \in C^{2,\alpha}(\mathcal{R}_T)
\]
is of class $C^1$ and $\f_{0,p,q}$ is the null function for every $(p, q) \in B_R$, we have that  
\begin{equation}\label{eq:phiprop}
\sup_{(p,q)\in B_R} \left\| \f_{\e,p,q}\right \|_{C^2(\mathcal{R}_T)} \leq C |\e|\,,
\end{equation}
where $C>0$ independent of $\e$, and 
\begin{equation}\label{eq:dephidpprop}
\sup_{(p,q)\in B_R} \left\| \frac{\partial \f_{\e,p,q}}{\partial p}\right \|_{C^2(\mathcal{R}_T)} \to 0\,, \quad \text{ as }\e \to 0\,.  
\end{equation}
Therefore 
\begin{equation}\label{eq:ex13}
\sup_{(p,q)\in B_R}\left\|\frac{\partial}{\partial p }(\tilde{H}_{\varepsilon}(U_{p,q} + \f_{\e,p,q}N_a)- \tilde{H}_{\varepsilon}(U_{p,q}) ) \right\|_{C^0(\mathcal{R}_T)} \to 0\,, \quad \text{ as }\e \to 0\,. 
\end{equation}

Now, setting
\[
\mathfrak{R_2}(\f_{\e, p, q}) := \mathfrak{M}(X_a + \f_{\e,p,q}N_a) - 1 - \mathfrak{L}_a \f_{\e,p,q}\,,
\]
we have
\begin{equation}\label{eq:derM}
\frac{\partial}{\partial p }\mathfrak{M}(X_a + \f_{\e,p,q}N_a) = \mathfrak{L}_a \frac{\partial \f_{\e,p,q}}{\partial p} + \frac{\partial}{\partial p} \mathfrak{R_2}(\f_{\e, p, q})\,.
\end{equation}
Since $\f_{\e,p,q} \in \X_T$ for any $\e,p,q$, then $ \frac{\partial \f_{\e,p,q}}{\partial p}(\pm T, \cdot) = 0$ on $\S^1$. Therefore, integrating by parts we find 
\begin{equation}\label{eq:kerdp}
 \int_{\mathcal{R}_T}x_a^2 \mathfrak{L}_a \left(\frac{\partial \f_{\e,p,q}}{\partial p}\right) w_{a,1}\cos \theta \, dt\, d \theta =  \int_{\mathcal{R}_T}x_a^2 \frac{\partial \f_{\e,p,q}}{\partial p}\,\mathfrak{L}_a ( w_{a,1} \cos \theta) \, dt\, d \theta = 0\,,
\end{equation}
where we used that $w_{a,1} \cos \theta \in \ker({\mathfrak{L}_{a})\big|}_{\X_T}$. Hence, from \eqref{eq:derM} and \eqref{eq:kerdp} we obtain
\begin{equation}\label{eq:ex14}
\int_{\mathcal{R}_T}x_a^2 \left( \frac{\partial}{\partial p }\mathfrak{M}(X_a + \f_{\e,p,q}N_a)\right) w_{a,1} \cos \theta \, dt\, d \theta = \int_{\mathcal{R}_T}x_a^2 \frac{\partial}{\partial p }\mathfrak{R_2}(\f_{\e,p,q}) w_{a,1} \cos \theta \, dt\, d \theta\,.
\end{equation}

A careful analysis of {\cite[Proposition B.1, B.2]{CIM}} shows that $\mathfrak{R_2}$ can be expressed, via Taylor expansion, as a linear combination of monomials of degree greater of equal than two with respect to $\f_{\e,p,q}$ and its derivatives with respect to $t, \theta$ up to the second order, having as coefficients smooth functions in $[-T, T]$ independent of $\e, p, q $. Therefore, differentiating with respect to $p$ we infer that
\[
\left\|\frac{\partial}{\partial p} \mathfrak{R_2}(\f_{\e, p, q})\right\|_{C^0(\mathcal{R}_T)} \leq C \left\| \frac{\partial \f_{\e,p,q}}{\partial p}\right \|_{C^2(\mathcal{R}_T)}\left\| \f_{\e,p,q}\right \|_{C^2(\mathcal{R}_T)}\,,
\]
where the constant $C>0$ is independent of $\e, p,q$. 
Then, using \eqref{eq:phiprop} and \eqref{eq:dephidpprop} we deduce that  
\begin{equation}\label{eq:ex15}
\e^{-1}\sup_{(p,q)\in B_R}\left\|\frac{\partial}{\partial p} \mathfrak{R_2}(\f_{\e, p, q})\right\|_{C^0(\mathcal{R}_T)} \to 0\,, \quad \text{ as }\e \to 0\,. 
\end{equation}

Putting together \eqref{eq:ex11}, \eqref{eq:ex12} and \eqref{eq:ex13}--\eqref{eq:ex15}, we find 
\[
\frac{\partial \F_1}{\partial p }(\e, p, q ) = \frac{\partial^2 M}{\partial p^2 }(p,q) + \mathfrak{R_3}(\e,p,q)\,,
\]
where $\mathfrak{R_3}$ is a continuous function such that 
\[
\sup_{(p,q)\in B_R}\mathfrak{R_3}({\e, p, q}) \to 0\,, \quad \text{ as }\e \to 0\,. 
\]

Using the same argument we can compute the whole Jacobian matrix of $\F$, obtaining that 
\[
J_{(p,q)}\F(\e, p, q) = D^2M(p, q) + \mathfrak{R}(\e, p, q), 
\]
where $D^2$ denotes the Hessian matrix, and the matrix $\mathfrak R = (\mathfrak{R}_{ij})$ is such that
\[
\sup_{(p,q)\in B_R}\mathfrak{R}_{i,j}(\e, p, q)\to 0 \quad \text{ as }\e \to 0\,. 
\] 

From this we conclude that $\F \in C^{1}([-\e, \e] \times B_R)$ and
\[
J_{(p,q)}\F(0, p_0, q_0) = D^2M(p_0, q_0)\,. 
\]
In fact, since $(p_0, q_0)$ is a nondegenerate critical point of $M$, then $J_{(p,q)}\F(0, p_0, q_0)$ is invertible. 

Applying the implicit function theorem to $\F$, we get that there exist $0 < \e_0 \leq \bar \e$ and uniquely determined $C^1$ functions
\[p:\left[
\begin{aligned}
&(-\e_0, \e_0) \to B_R\\
&\e \mapsto p_\e\,,
\end{aligned}\right.
\qquad 
q: \left[
\begin{aligned}
&(-\e_0, \e_0)\to B_R\\
&\e \mapsto q_\e\,,
\end{aligned}\right.
\] 
such that 
\[
p(0)  = p_0, \quad q(0) = q_0, \quad \F(\e, p_\e, q_\e) = 0\, \quad \text{ for any }\e \in (-\e_0, \e_0)\,.
\]
In particular, by definition of $\mathcal F$, it holds that 
\[
\lambda_{\e , 1}(p_\e, q_\e) = \lambda_{\e,2}(p_\e, q_\e) = 0 \quad \text{ for any }|\e|\leq \e_0\,, 
\]
and by \eqref{eq:c2ex1} we conclude that 
\[
\mathfrak{M}(U_{p_\e,q_\e} + \f_{\e,p_\e,q_\e}N_a) = H_{\varepsilon}(U_{p_\e,q_\e}+ \f_{\e,p_\e,q_\e}N_a)\,\quad \text{ for any }|\e|\leq \e_0\,.
\]
Finally, setting $\f_\e := \f_{\e,p_\e,q_\e}$, we have that the triple $(p_\e, q_\e, \f_\e)$ solves \eqref{M=H} and Theorem \ref{T:existence2}, $(ii)$ is proved.  \qed

\begin{appendices}

\section{\!\!\!\!\!\!: Jacobi elliptic functions}\label{A:Jacobi}
In this section we present some definitions and results regarding the Jacobi elliptic functions. For more details, we refer to \cite[Chapters 15,16]{AS}. 

The {\it incomplete elliptic integrals} of the first and second kind of parameter $m$ are defined, respectively, as 
\begin{equation}\label{eq:incellint}
F(s |m) := \int_0^{s}\frac{1}{\sqrt{1-m\sin^2 \theta}}\,d\theta \quad \quad E(s|m) :=\int_{0}^{s}\sqrt{1-m\sin^2\theta}\, d\theta\,, \qquad s \in \R\,, m \in [0,1)\,.
\end{equation}

The {\it Jacobi Amplitude} function of parameter $m$, denoted by $\mathrm{am}(s|m)$, is defined as the inverse function of $F$ with respect to $s$, namely
\begin{equation}\label{eq:am}
s = F(\am(s|m)|m)\,,\qquad s\in \R\,,\ m \in [0,1).
\end{equation}

The {\it Delta Amplitude} function of parameter $m$ is defined as 
\[
\mathrm{dn}(s|m) := \frac{\partial}{\partial s}\mathrm{am}(s|m) \,,\qquad s\in \R\,,\ m \in [0,1)\,.
\]
For convenience, we also recall the following alternative representation
\begin{equation}\label{eq:dn}
\mathrm{dn}(s|m) =\sqrt{1-m\sin^2 (\mathrm{am}(s|m))}\,,\qquad s\in \R\,,\ m \in [0,1).
\end{equation}

Finally, the {\it complete elliptic integrals} of the first and second kind are defined, respectively, as 
\[
K(m) := F\left(\frac{\pi}{2}\Big|m\right) \quad \quad E(m) := E\left(\frac{\pi}{2}\Big|m\right)\,, \qquad m \in [0,1)\,.
\]

Here we list some known properties of the Jacobi elliptic functions, see also \cite[22.13, 22.16]{DLMF}.
\begin{gather}
\mathrm{am}(0|m) = 0\,, \qquad \mathrm{am}(K(m)|m) = \frac{\pi}{2}\,,\label{eq:extrema}\\
\frac{\partial}{\partial s}\mathrm{dn}(s|m) = -m \sin(\mathrm{am}(s|m))\cos(\am(s|m))\,.\label{eq:dnds}\\
\label{eq:x2int}
\int_0^{s}\dn(\tau|m)^2\,d\tau = E(\am(s|m)|m)\,,\\
\label{eq:1sux2int}
\int_0^s \frac{d\tau}{\dn(\tau|m)^2} = \frac{1}{1-m}\left[E(\am(s|m)|m)+\frac{\frac{\partial}{\partial s}\mathrm{dn}(s|m) }{\dn(s|m)}\right]\,.
\end{gather}
Notice that from \eqref{eq:extrema}, \eqref{eq:dnds} and \eqref{eq:1sux2int} it follows that
\begin{equation}\label{eq:dnintcei}
\int_0^{K(m)}\dn(\tau|m)^2\,d\tau = E(m)\,, \qquad \qquad\int_0^{K(m)} \frac{d\tau}{\dn(\tau|m)^2} = \frac{E(m)}{1-m}\,.
\end{equation}

In the next Lemma we compute the derivative of the Delta amplitude function with respect to $m$.

\begin{Lemma} For any $s \in \R$, $m \in [0,1)$, it holds
\begin{equation}
\begin{aligned}
\label{eq:dndk}
\frac{\partial }{\partial m}\mathrm{dn}(s|m) &= \frac{1}{2m}\left[\mathrm{dn}(s|m) - \frac{1}{\mathrm{dn}(s|m)}\right]\\
&+\frac{\frac{\partial }{\partial s}\mathrm{dn}(s|m)}{2m}\left[ s - \frac{\frac{\partial}{\partial s}\mathrm{dn}(s|m)}{(1-m)\mathrm{dn}(s|m)} - \frac{1}{1-m}\int_0^s\mathrm{dn}(\tau|m)^2 d\tau \right]\,.
\end{aligned}
\end{equation}
\end{Lemma} 
\begin{Proof} Differentiating \eqref{eq:incellint} with respect to $m$ and exploiting the following relations
\begin{equation}\label{eq:thetarel} 
\begin{gathered}
\frac{\sin^2 \theta}{\sqrt{1-m\sin^2 \theta}} =-\frac{1}{m}\left[\sqrt{1-m\sin^2\theta} - \frac{1}{\sqrt{1-m\sin^2\theta}}\right]\,,\\
\frac{1}{(1-m\sin^2\theta)^{3/2}} =  \frac{\sqrt{1-m\sin^2\theta}}{1-m} -\frac{m}{1-m}\frac{\partial}{\partial \theta}\left(\frac{\sin \theta \cos \theta}{\sqrt{1-m\sin^2\theta}} \right)\,,
\end{gathered}
\end{equation}
we get that 
\begin{equation}\label{eq:dFdm}
\frac{\partial }{\partial m} F(s|m) = -\frac{1}{2m}\left[F(s|m) + \frac{\sin(s) \cos(s) }{(1-m)\sqrt{1-m\sin^2(s)}}-\frac{E(s|m)}{1-m}\right]\,.
\end{equation}

Next, from \eqref{eq:am}, \eqref{eq:dnds}, \eqref{eq:x2int} and \eqref{eq:dFdm} we find 
\begin{equation}
\label{eq:amdk}
\frac{\partial}{\partial m}\am(s|m) = \frac{\mathrm{dn}(s|m)}{2m}\left[ s -\frac{\frac{\partial}{\partial s}\mathrm{dn}(s|m)}{(1-m)\mathrm{dn}(s|m)}-\frac{1}{1-m} \int_{0}^{s}\mathrm{dn}(\tau|m)^2 d\tau  \right]\,.
\end{equation}

Finally, differentiating \eqref{eq:dn} with respect to $m$ and using \eqref{eq:dnds} we obtain
\begin{equation}\label{eq:ddmdn}
\frac{\partial}{\partial m}\left[ \mathrm{dn}(s|m)\right] =- \frac{\sin^2 (\am(s|m))}{2\,\mathrm{dn}(s|m)} + \frac{\frac{\partial}{\partial s}\mathrm{dn}(s|m)}{\mathrm{dn}(s|m)}\frac{\partial}{\partial m}\am(s|m) \,.
\end{equation}
Combining \eqref{eq:thetarel}, \eqref{eq:amdk} and \eqref{eq:ddmdn} we deduce \eqref{eq:dndk}. The Lemma is proved. \qed
\end{Proof}

\section{\!\!\!\!\!\!: Volume functional}\label{A:areavol}
Let $\mathcal{S}_{T}$ be the class of maps $X=X(t,\theta)\colon[-T,T]\times\R/2\pi\Z\to\R^{3}$ of class $C^{2}$, $2\pi$-periodic with respect to $\theta$, with 
\[
X_{t}(t,\theta)\wedge X_{\theta}(t,\theta)\ne 0\quad\forall(t,\theta)\in[-T,T]\times\R
\]
and such that $X\cdot\mathbf{e}_{1}$ and $X\cdot\mathbf{e}_{2}$ are even with respect to $t$, whereas $X\cdot\mathbf{e}_{3}$ is odd with respect to $t$.

We observe that for every $a\in(-1,0)\cup(0,\infty)$ the mapping $X_{a}|_{[-T,T]\times\R}$ belongs to $\mathcal{S}_{T}$. Moreover, when $\varphi=\varphi(t,\theta)\colon[-T,T]\times\R/2\pi\Z\to\R^{3}$ is a function of class $C^{2}$, $2\pi$-periodic with respect to $\theta$, even with respect to $t$, and with $\|\varphi\|_{C^{0}(\mathcal{R}_T)}$ small enough, then $X_{a}+\varphi N_{a}\in\mathcal{S}_{T}$, as it follows by definition and by \eqref{eq:extcond}. 

\medskip

Let $H\colon\R^{3}\to\R$ be a continuous function and let $Q\colon\R^{3}\to\R^{3}$ be any vector field of class $C^{1}$ and such that $\mathrm{div}\,Q=H$ in $\R^{3}$. We point out that such a vector field $Q$ always exists (see for instance \eqref{eq:Q}) but, in general, it is not uniquely determined.    
For any choice of $Q$ we set
\begin{equation}\label{H-vol-def}
\mathcal{V}_{T, Q}(X)=\int_{\mathcal{R}_{T}}Q(X)\cdot X_{t}\wedge X_{\theta}\,dt\,d\theta\quad\forall X\in\mathcal{S}_{T}\,.
\end{equation}
\begin{Remark}
We observe that when $X = X_a$ with $a\in(-1,0)$, the mapping $X_{a}$ is a parameterization of an unduloid whose axis of revolution is the $z$-axis. Let 
\[
\Theta_{a,T}:=\{(\rho x_{a}(t)\cos\theta,\rho x_{a}(t)\sin\theta,z_{a}(t))\in\R^{3}\mid \rho\in[0,1]\,,~|t|\le T, \theta \in [-\pi, \pi]\}
\]
be the region of $\R^{3}$ enclosed by the unduloid and the horizontal planes $z=\pm z_{a}(T)$. 
Then the $H$-weighted volume of $\Theta_{a,T}$ is given by
\begin{equation}\label{H-vol-Xa}
\mathrm{vol}_{H}(\Theta_{a,T})=-\mathcal{V}_{T, Q}(X_{a})+\int_{\Sigma^{+}_{a,T}}Q_{3}\,d\sigma-\int_{\Sigma^{-}_{a,T}}Q_{3}\,d\sigma
\end{equation}
where $\Sigma_{a,T}^{\pm}:=\{(x,y,z)\in\R^{3}\mid x^{2}+y^{2}\le x_{a}^{2}(T)\,,~z=z_{a}(\pm T)\}$. Indeed, by the divergence theorem we have
\[
\mathrm{vol}_{H}(\Theta_{a,T})=\int_{\Theta_{a,T}}H(x,y,z)\,dx\,dy\,dz=\int_{\partial \Theta_{a,T}}Q\cdot\nu\,d\sigma
\]
where $Q:\R^3 \to \R^3$ is any vector field such that $\mathrm{div}\,Q=H$ in $\R^{3}$ and $\nu$ is the outer normal to $\partial \Theta_{a,T}$. Notice that $\partial \Theta_{a,T}=\Sigma_{a,T}\cup\Sigma_{a,T}^{+}\cup\Sigma_{a,T}^{-}$, where $\Sigma_{a,T}=X_{a}(\mathcal{R}_{T})$. In addition $\nu=-N_{a}$ on $\Sigma_{a,T}$ whereas $\nu=\pm\mathbf{e}_{3}$ on $\Sigma_{a,T}^{\pm}$. Hence
\[
\mathrm{vol}_{H}(\Theta_{a,T})=-\int_{\mathcal{R}_{T}}Q(X_{a})\cdot N_{a}|(X_{a})_{t}\wedge(X_{a})_{\theta}|\,dt\,d\theta+\int_{\Sigma^{+}_{a,T}}Q_{3}\,d\sigma-\int_{\Sigma^{-}_{a,T}}Q_{3}\,d\sigma
\]
and \eqref{H-vol-Xa} follows from \eqref{H-vol-def} and from the definition of $N_{a}$.
\end{Remark}

Fixing $X\in\mathcal{S}_{T}$, we call a \emph{variation of $X$} a curve $s\mapsto X(s)\in\mathcal{S}_{T}$ of class $C^{2}$ defined in an open interval $(-s_{0},s_{0})\subset\R$ and  such that $X(0)=X$. For convenience, we set $\dot{X}:=\left.\frac{\partial X(s)}{\partial s}\right|_{s=0}$, and we observe that $\dot X$ shares the same symmetries of $X$.

\begin{Lemma}
\label{L:first-variation-volume}
Let $X\in\mathcal{S}_{T}$ and let $s\mapsto X(s)\in\mathcal{S}_{T}$ be a variation of $X$. Then
\[
\begin{aligned}
\left.\frac{\partial}{\partial s}\left[\mathcal{V}_{T, Q}(X(s))\right]\right|_{s=0}&=\int_{\mathcal{R}_{T}}H(X)\dot{X}\cdot X_{t}\wedge X_{\theta}\,dt\,d\theta\\
&-\int_{-\pi}^{\pi}\![(Q(X)\cdot X_{\theta}\wedge \dot{X})(T,\theta)-(Q(X)\cdot X_{\theta}\wedge \dot{X})(-T,\theta)]\,d\theta\,.
\end{aligned}
\]
\end{Lemma}
\begin{Proof} 
By direct computation we have that 
\[\begin{split}
\left.\frac{\partial}{\partial s}\left[\mathcal{V}_{T, Q}(X(s))\right]\right|_{s=0}&=\int_{\mathcal{R}_{T}}[JQ(X)\dot{X}]\cdot X_{t}\wedge X_{\theta}\,dt\,d\theta\\
&+\int_{\mathcal{R}_{T}}Q(X)\cdot\dot{X}_{t}\wedge X_{\theta}\,dt\,d\theta+\int_{\mathcal{R}_{T}}Q(X)\cdot X_{t}\wedge \dot{X}_{\theta}\,dt\,d\theta\,,
\end{split}\]
where $JQ$ is the Jacobian matrix of $Q$.  
Then, applying the algebraic identity
\[ (M\mathbf{a})\cdot\mathbf{b}\wedge\mathbf{c}+
\mathbf{a}\cdot(M\mathbf{b})\wedge\mathbf{c}+
\mathbf{a}\cdot\mathbf{b}\wedge(M\mathbf{c})=
(\mathrm{tr}\,M)\mathbf{a}\cdot\mathbf{b}\wedge\mathbf{c}\,, \quad \forall\, \mathbf{a}, \mathbf{b}, \mathbf{c} \in \R^3\,, \forall\, M \in \R^{3\times 3}
\]
with $\mathbf{a}=\dot X$, $\mathbf{b}=X_{t},\mathbf{c}=X_{\theta}$, $M=JQ(X)$, we obtain 
\[\begin{split}
\left.\frac{\partial}{\partial s}\left[\mathcal{V}_{T,Q}(X(s))\right]\right|_{s=0}=
&\int_{\mathcal{R}_{T}}H(X)\dot{X}\cdot X_{t}\wedge X_{\theta}\,dt\,d\theta-\int_{\mathcal{R}_{T}}\dot{X}\cdot\left[X_{t}\wedge(Q(X))_{\theta}+(Q(X))_{t}\wedge X_{\theta}\right]\,dt\,d\theta\\
&+\int_{\mathcal{R}_{T}}Q(X)\cdot\left[\dot{X}_{t}\wedge X_{\theta}+X_{t}\wedge \dot{X}_{\theta}\right]\,dt\,d\theta\,.
\end{split}\]
To conclude it is sufficient to use the identity
\begin{equation}\label{Wente}
\int_{\mathcal{R}_{T}}Y\cdot[X_{t}\wedge Z_{\theta}+Z_{t}\wedge X_{\theta}]\,dt\,d\theta=\int_{\mathcal{R}_{T}}Z\cdot[Y_{t}\wedge X_{\theta}+X_{t}\wedge Y_{\theta}]-\int_{-\pi}^{\pi}\Big[Z\cdot X_{\theta}\wedge Y\Big]_{t=-T}^{t=T}d\theta
\end{equation}
with $Y=\dot X$ and $Z=Q(X)$. It remains to prove \eqref{Wente}. For this, integrating by parts we have 
\[\begin{split}
\int_{\mathcal{R}_{T}}Y\cdot[X_{t}\wedge Z_{\theta}+Z_{t}\wedge X_{\theta}]\,dt\,d\theta&=
\int_{\mathcal{R}_{T}}Z_{\theta}\cdot Y\wedge X_{t}\,dt\,d\theta+\int_{\mathcal{R}_{T}}Z_{t}\cdot X_{\theta}\wedge Y\,dt\,d\theta\\
&=
\int_{-T}^{T}\Big[Z\cdot Y\wedge X_{t}\Big]_{\theta=-\pi}^{\theta=\pi}dt-\int_{\mathcal{R}_{T}}Z\cdot[Y_{\theta}\wedge X_{t}+Y\wedge X_{t\theta}]\,dt\,d\theta\\
&\quad
+\int_{-\pi}^{\pi}\Big[Z\cdot X_{\theta}\wedge Y\Big]_{t=-T}^{t=T}d\theta-\int_{\mathcal{R}_{T}}Z\cdot[X_{\theta t}\wedge Y+X_{\theta}\wedge Y_{t}]\,dt\,d\theta\,.
\end{split}\] 
Then, since $X$, $Y$ and $Z$ are $2\pi$-periodic with respect to $\theta$, it follows that $\big[Z\cdot Y\wedge X_{t}\big]_{\theta=-\pi}^{\theta=\pi}=0$ and \eqref{Wente} readily follows. The proof of the Lemma is complete. \qed

\end{Proof}

\end{appendices}

\end{document}